\documentclass[12pt]{article}
\usepackage{amscd}
\usepackage{latexsym}
\usepackage{mathrsfs}
\usepackage{pst-all,multido,ifthen}
\usepackage{amsmath}
\usepackage{amssymb}
\usepackage[all]{xy}
\usepackage[english]{babel}
\usepackage{mathtools}
\usepackage{upgreek}
\usepackage{scalerel}
\usepackage{stackengine,wasysym}

\stackMath
\usepackage{tikz}

\newcommand\reallywidehat[1]{%
\savestack{\tmpbox}{\stretchto{%
  \scaleto{%
    \scalerel*[\widthof{\ensuremath{#1}}]{\kern-.6pt\bigwedge\kern-.6pt}%
    {\rule[-\textheight/2]{1ex}{\textheight}}%WIDTH-LIMITED BIG WEDGE
  }{\textheight}% 
}{0.5ex}}%
\stackon[1pt]{#1}{\tmpbox}%
}

\DeclarePairedDelimiter{\ceil}{\lceil}{\rceil}
\newtheorem{thm}{Theorem}

\newtheorem{prop}{Proposition}
\newtheorem{cor}{Corollary}
\newtheorem{lemma}{Lemma}

\newtheorem{fact}{Fact}
\newtheorem{claim}{Claim}
\newtheorem{ex}{Example}
\def\endproof{$\hfill \square$}
\def\Spec{\textup{Spec}}
\def\BT{\textup{BT}}
\def\Max{\textup{Max}}
\def\depth{\textup{depth}}
\def\Ker{\textup{Ker}}
\def\codim{\textup{codim}}
\def\Frac{\textup{Frac}}
\def\red{\textup{red}}
\def\dim{\textup{dim}}

\begin{document}
\title{Purity for Barsotti--Tate groups in some mixed characteristic situations}
\author{Ofer Gabber and Adrian Vasiu}
\maketitle

\noindent
{\bf ABSTRACT.} Let $p$ be a prime. Let $R$ be a regular local ring of dimension $d\ge 2$ whose completion is isomorphic to $C(k)[[x_1,\ldots,x_d]]/(h)$, with $C(k)$ a Cohen ring with the same residue field $k$ as $R$ and with $h\in C(k)[[x_1,\ldots,x_d]]$ such that its reduction modulo $p$ does not belong to the ideal $(x_1^p,\ldots,x_d^p)+(x_1,\ldots,x_d)^{2p-2}$ of $k[[x_1,\ldots,x_d]]$. We extend a result of Vasiu--Zink (for $d=2$) to show that each Barsotti--Tate group over $\Frac(R)$ which extends to every local ring of $\Spec(R)$ of dimension $1$, extends uniquely to a Barsotti--Tate group over $R$. This result corrects in many cases several errors in the literature. As an application, we get that if $Y$ is a regular integral scheme such that the completion of each local ring of $Y$ of residue characteristic $p$ is a formal power series ring over some complete discrete valuation ring of absolute ramification index $e\le p-1$, then each Barsotti--Tate group over the generic point of $Y$ which extends to every local ring of $Y$ of dimension $1$, extends uniquely to a Barsotti--Tate group over $Y$. 



\section{Introduction}\label{S1}
Let $R$ be a local noetherian ring with residue field $k$. Let $X:=\Spec(R)$ and let $U:=X\setminus\Spec(k)$ be the punctured spectrum of $R$. Let $p$ be a prime number. The notation $R$, $k$, $X$, $U$ and $p$ is fixed throughout this article. 

We will abbreviate a Barsotti--Tate group (i.e., a $p$-divisible group) by $\BT$ and a truncated Barsotti--Tate group of level $n\in\mathbb N^*$ by $\BT_n$. We recall from \cite{VZ2}, Def. 2 that when $\depth(R)\ge 2$, $R$ is called {\it $p$-quasi-healthy} if each $\BT$ over $U$ extends to a $\BT$ over $X$; from \cite{G2}, Exp. III, Cor. 3.5 applied to coherent sheaves defined by the structure sheaves of $\BT_n$'s over $X$ we get that such an extension is unique up to unique isomorphism. 

We recall that the first examples of $\BT$'s over $X$ are obtained by considering an abelian scheme $\mathcal A$ over $X$ and by taking a direct summand of the $\BT$ $\mathcal A[p^{\infty}]$ of $\mathcal A$, i.e., of the inductive system $\mathcal A[p^n]=\Ker(p^n:\mathcal A\rightarrow\mathcal A)$ indexed by $n\in\mathbb N^*$. If $R$ is complete and $k$ is an algebraically closed field of characteristic $p$, then the converse holds, i.e., each $\BT$ over $X$ is a direct summand of $\mathcal A[p^{\infty}]$ for some abelian scheme $\mathcal A$ over $X$ (this follows from \cite{V3}, Prop. 5.3.3).

If $p=0$ in $R$ and $\depth(R)\ge 2$, then one can check that $R$ is not $p$-quasi-healthy using a variant of Moret-Bailly's example in \cite{FC}, p. 192.\footnote{We consider a homomorphisms $\pmb{\alpha}_{{\bf p},X}\rightarrow\pmb{\alpha}_{{\bf p},X}^{\dim(R)}$ defined by a system of parameters of $R$ (it is a closed embedding over $U$ but not over $X$) and (as in a theorem of Raynaud) an embedding of $\pmb{\alpha}_{{\bf p},X}^{\dim(R)}$ as a closed subgroup scheme of an abelian scheme $\mathcal A$ over $X$. If $\mathcal D_U:=\mathcal A_U[p^{\infty}]/\pmb{\alpha}_{{\bf p},U}$ extends to a $\BT$ $\mathcal D$ over $X$, then by the depth assumption on $R$ the isogeny $\mathcal A_U[p^{\infty}]\rightarrow\mathcal D_U$ extends to a morphism $\mathcal A[p^{\infty}]\rightarrow\mathcal D$, which must be an isogeny (cf. \cite{CCO}, Prop. 3.3.8 and the last part of Ex. 3.3.10), so its kernel is a finite flat subgroup scheme of $\pmb{\alpha}_{{\bf p},X}^{\dim(R)}$, which leads to a contradiction.} %Raynaud’s theorem of [BBM, 3.1.1]

In \cite{VZ2}, Thm. 3 it is shown that there are large classes of $p$-quasi-healthy regular rings of dimension $2$ and mixed characteristic $(0,p)$. This result already has applications to the uniqueness of integral canonical models of Shimura varieties and to the existence of new classes of N\'eron models (see \cite{VZ2}, Cor. 30 and Thm. 31), to the study of moduli spaces of polarized K3 surfaces (see \cite{MP1}, \cite{MP2}), to the study of crystalline representations associated to fundamental groups of suitable formally smooth algebras over complete discrete valuation rings of absolute ramification index at most $p-2$ (see \cite{LM}, \cite{MoYS}), and to the computation of certain \'etale cohomology classes (see \cite{GV}). 

But no examples of $p$-quasi-healthy regular rings of mixed characteristic $(0,p)$ are (correctly) proved to exist in the literature for dimensions at least $3$; this is so due to the fact that all claims in the literature for dimensions at least $3$ rely on the erroneous argument in \cite{FC}, Ch. V, Sect. 6, end of p. 183 and top of p. 184. The difficulty encountered consists in being able to provide examples of dimension $3$ as the Grothendieck--Messing deformation theory for $\BT$'s allows a relatively direct passage from dimension $3$ to higher dimensions (see Lemma \ref{L5}). 

Our goal is to obtain purity results for $\BT$'s over regular schemes. More precisely, we will provide many examples of $p$-quasi-healthy regular local rings of dimension at least $3$ and we will use them to get as well new examples of faithfully flat $\Spec(\mathbb Z_{(p)})$-schemes which are $p$-healthy regular in the sense of either \cite{V1}, Def. 3.2.1 9) or \cite{VZ2}, Def. 1. Our main result is the following theorem which extends \cite{VZ2}, Thm. 3. 

\begin{thm}\label{T1}
Let $R$ be a regular local ring of dimension $d\ge 1$ and mixed characteristic $(0,p)$ which satisfies the following condition:

\medskip
{\bf ($\natural$)} The completion of $R$ is isomorphic to $C(k)[[x_1,\ldots,x_d]]/(h)$, where $C(k)$ is a Cohen ring of $k$ and where $h\in C(k)[[x_1,\ldots,x_d]]$ is such that its reduction $\bar h$ modulo $p$ does not belong to the ideal $(x_1^p,\ldots,x_d^p)+(x_1,\ldots,x_d)^{2p-2}$ of $k[[x_1,\ldots,x_d]]$. 

\medskip
Then the following two properties hold:

\medskip
{\bf (a)} If a $\BT$ over $\Frac(R)$ extends to each local ring of $X=\Spec(R)$ of dimension $1$, then it extends uniquely (up to unique isomorphism) to a $\BT$ over $X$. 

\smallskip
{\bf (b)} If $d\ge 2$, then the regular local ring $R$ is $p$-quasi-healthy. 
\end{thm}

Clearly (a) implies (b) and the fact that $X$ is also $p$-healthy regular. Condition ($\natural$) is stable under generization, cf. Proposition \ref{P5}. Based on this and on the classical purity theorem of Zariski, Nagata and Grothendieck (see \cite{G2}, Exp. X, Thm. 3.4 (i)]), we get that in fact (a) and (b) are equivalent (see Section \ref{S7}). 

For $d=2$, Theorem \ref{T1} (b) is proved in \cite{VZ2}, Thm. 3. Theorem \ref{T1} (b) and \cite{VZ2}, Prop. 23 (b) imply \cite{VZ2}, Thm. 3 and that $R$ for $d\ge 2$ is also a quasi-healthy regular ring in the sense of \cite{VZ2}, Def. 2, which means that each abelian scheme over $U$ extends (uniquely up to a unique isomorphism) to an abelian scheme over $X$.\footnote{It is an open problem when this extension property holds for polarized K3 surfaces.} In this way we get precisely those examples of quasi-healthy regular rings one gets based only on \cite{VZ2}, Thm. 3 and Lem. 24. This is so as based on Propositions \ref{P3} and \ref{P4} we easily get that for $d\ge 2$ the condition on $R$ in \cite{VZ2}, Thm. 3 is equivalent to ($\natural$).

For $d\ge 3$, we prove Theorem \ref{T1} (b) by induction on $d\ge 3$ (see Section \ref{S6}) and the hard part is the case when $d=3$ (see Subsection \ref{SS61}). The proof of Theorem \ref{T1} (b) for $d=3$ involves a study of the blow up of $X$ along its closed point $\Spec(k)$ and relies heavily on several important results and ideas. 

Firstly, the proof of Theorem \ref{T1} (b) for $d=3$ relies on an application (see Proposition \ref{P2}) of Raynaud's complement to Tate's extension theorem of \cite{T}, Thm. 4 obtained in \cite{R}, Prop. 2.3.1 and reproved in Proposition \ref{P1} based on a refinement of \cite{VZ3}, Cor. 4 proved in Lemma \ref{L1}. Secondly, the proof of Theorem \ref{T1} (b) for $d=3$ relies on the case $d=2$ of Theorem \ref{T1} (b) proved in \cite{VZ2}. Thirdly, it relies on the particular case $(N,l^{\prime},\mathcal C)=(2,l,\mathbb P^1_l)$ of the general theorem below which we think is of interest in its own. 

\begin{thm}\label{T2}
Let $l$ be a field of characteristic $p$. Let $N\in\mathbb N^*$ be an integer. For each $n\in\mathbb N^*$ let $\mathcal D_n$ be a $\BT_n$ over an open subscheme $\mathcal U_n$ of $\mathbb P^N_l$. We assume that we have a chain of inclusions $\mathcal U_1\supset\mathcal U_2\supset\cdots\supset\mathcal U_n\supset\cdots$ and that for all $n\in \mathbb N^*$ we have an identification $\mathcal D_{n+1}[p^n]=\mathcal D_{n,\mathcal U_{n+1}}$. We also assume that there exists a field extension $l^{\prime}$ of $l$ and a closed subscheme $\mathcal C$ of $\mathbb P^N_{l^{\prime}}$ of dimension greater than $0$ and contained in $\mathcal U_{n,l^{\prime}}$ for all $n\in\mathbb N^*$ such that the inductive system $\mathcal C\times_{\mathcal U_n} \mathcal D_n$ indexed by $n\in\mathbb N^*$ is a constant $\BT$ over $\mathcal C$, i.e., it is isomorphic to the pullback to $\mathcal C$ of a $\BT$ over $\Spec(l^{\prime})$. Then for each $n\in\mathbb N^*$, $\mathcal D_n$ extends uniquely (up to unique isomorphism) to a constant $\BT_n$ $\mathcal D_n^+$ over $\mathbb P^N_l$ isomorphic to $G[p^n]_{\mathbb P^N_l}$ for a suitable $\BT$ $G$ over $\Spec(l)$ and the identification $\mathcal D_{n+1}[p^n]=\mathcal D_{n,\mathcal U_{n+1}}$ extends to an identification $\mathcal D_{n+1}^+[p^n]=\mathcal D_n^+$ which is the pullback to $\mathbb P^N_l$ of the identification $G[p^{n+1}][p^n]=G[p^n]$ over $\Spec(l)$ (therefore the $\BT$ over the stable under generization, pro-constructible subset $\mathcal U_{\infty}:=\cap_{n=1}^{\infty} \mathcal U_n$ of $\mathbb P^N_l$ induced naturally by the $\mathcal D_n$'s, extends uniquely up to unique isomorphism to a constant $\BT$ over $\mathbb P^N_l$ isomorphic to $G_{\mathbb P^N_l}$\footnote{Note that $\mathcal U_{\infty}$ is a locally ringed space, and one can define a $\BT$ over a locally ringed space $(\mathcal Y,\mathcal O_{\mathcal Y})$ as a projective system of commutative and cocommutative Hopf $\mathcal O_{\mathcal Y}$-algebras $\mathcal H_n$ indexed by $n\in\mathbb N^*$, locally free of finite rank as 
$\mathcal O_{\mathcal Y}$-modules, such that for every point $z\in\mathcal Y$ the inductive system $\Spec((\mathcal H_n)_z)$ of finite flat group schemes indexed by $n\in\mathbb N^*$ and defined by the stalks at $z$ constitutes a $\BT$ over the local scheme $\Spec(\mathcal O_{\mathcal Y,z})$.}).
\end{thm}

The proof of Theorem \ref{T2} is presented in Section \ref{S5} and it relies on properties of rings of formal-rational functions established in \cite{HM}. Theorem \ref{T1} is proved in Sections \ref{S6} and \ref{S7}. Section \ref{S2} presents a few basic extension results that are required in the proof and the applications of Theorem \ref{T1}, including a descent lemma (see Lemma \ref{L3}) for $p$-quasi-healthy regular rings and the proof of Raynaud's complement to it. Section \ref{S3} groups together a few basic properties related to condition ($\natural$). Section \ref{S4} contains an elementary formal algebraic geometry property which is a standard application of the results of \cite{HM} and which plays a key role in the proof of Theorem \ref{T2}. 
 
The following corollary of Theorem \ref{T1} and its proof and of Bondarko's boundedness results for truncated $\BT$'s over discrete valuation rings of mixed characteristic $(0,p)$ (see \cite{B} and \cite{VZ3}) extends \cite{VZ2}, Cor. 5 and it can be used to correct the argument used in \cite{FC}, Ch. V, Sect. 6, end of p. 183 and top of p. 184 in many situations.

\begin{cor}\label{C1}
Let $Y$ be a regular integral scheme flat over $\Spec(\mathbb Z)$. We assume that for each point $z\in Y$ of characteristic $p$, condition ($\natural$) holds for the local ring $\mathcal O_{Y,z}$ of $z$ in $Y$; this holds if the completion of $\mathcal O_{Y,z}$ is formally smooth (e.g., is a formal power series ring) over a complete discrete valuation ring of absolute ramification index $e\in\{1,\ldots,p-1\}$. Then the following two properties hold:

\medskip
{\bf (a)} Each $\BT$ over the generic point of $Y$ which extends to every local ring of $Y$ of dimension $1$, extends uniquely (up to unique isomorphism) to a $\BT$ over $Y$.

\smallskip
{\bf (b)} Each $\BT$ over an open subscheme of $Y$ that contains $Y[\frac{1}{p}]$ and all points of $Y_{\mathbb F_p}$ of codimension $0$ in $Y_{\mathbb F_p}$ extends to a $\BT$ over $Y$ (thus if $Y$ is a faithfully flat $\Spec(\mathbb Z_{(p)})$-scheme, then it is $p$-healthy regular).
\end{cor}

The case $e<p-1$ of Corollary \ref{C1} (b) is a special case of \cite{FC}, Ch. V, Thm. 6.4' (which is incorrect already in dimension $2$) and was claimed more recently in \cite{Mo}, Subsect. 3.6.1 and in \cite{V1}, Rm. 3.2.2 3) and the last two lines of Subsubsect. 3.2.17. The case $e=1$ of Corollary \ref{C1} (b) was also claimed in \cite{V2}, Thm. 1.3. Unfortunately, the inductive passages from dimension $2$ to dimensions $3$ in all these four references relied on \cite{FC}, Ch. V, Sect. 6, argument on pp. 183--184, which, in the language of this paper, would in particular imply that for each integer $n\ge 2$ all local noetherian rings $C(k)[[x_1,x_2]]/(x_2^n)$ are $p$-quasi-healthy. The counterexample of \cite{VZ2}, Subsect. 5.1 shows that the last statement is false. 

The hypothesis of Corollary \ref{C1} (a) holds for each $\BT$ over the generic point of $Y$ whose truncations extend to finite flat group schemes over $Y$ (cf. Raynaud's extension theorem). Thus Corollary \ref{C1} (a) reobtains, refines and generalizes the first part of \cite{MoYS}, Thm. 1.2 which worked in certain formally smooth contexts over a discrete valuation ring of absolute ramification index $e<p-1$.

If $Y$ is as in Corollary \ref{C1} and faithfully flat over $\Spec(\mathbb Z_{(p)})$, then $Y$ is also a healthy regular scheme in the sense of \cite{V1}, Def. 3.2.1 2), i.e., each abelian scheme over an open subscheme of $Y$ that contains $Y_{\mathbb Q}$ and all points of $Y_{\mathbb F_p}$ of codimension $0$ in $Y_{\mathbb F_p}$ extends to an abelian scheme over $Y$ (cf. \cite{V2}, Prop. 4.1). From this and \cite{VZ2}, Lem. 29 one gets a second proof of the uniqueness of integral canonical models of Shimura varieties defined in \cite{V1}, Def. 3.2.3 6) over discrete valuation rings of mixed characteristic $(0,p)$ and absolute ramification index at most $p-1$ which was first proved in \cite{VZ2}, Cor. 30. 

The examples of \cite{VZ2}, Thm. 28 (i) and (ii) for $d=2$ can be easily adapted to provide many examples of regular local rings $R$ of arbitrary dimension $d\ge 2$ which are of mixed characteristic $(0,p)$ and are not $p$-quasi-healthy. For example, the regular complete local rings $C(k)[[x_1,\ldots,x_d]]/(p-ax_1^p)$ and $C(k)[[x_1,\ldots,x_d]]/(p-a\prod_{i=1}^d x_i^{p-1})$ are not $p$-quasi-healthy, where $a\in C(k)[[x_1,\ldots,x_d]]$; in particular, formal power series rings in $d-1\ge 1$ variables over a complete discrete valuation ring of absolute ramification index at least $p$ are not $p$-quasi-healthy. See \cite{GV} for details, generalizations, and the fact that if $R$ is regular henselian of dimension $2$, then condition ($\natural$) holds if and only if $R$ is $p$-quasi-healthy.

Complements to Corollary \ref{C1} and to Lemma \ref{L8} used in its proof are included in Section \ref{S8}. 

Let $\mathcal O_{\flat}$ be the structure sheaf of a scheme $\flat$. For each local noetherian ring $R$, let $\widehat{R}$ be its completion. If $R$ is an integral domain, let $K:=\Frac(R)$; so if $R$ is a discrete valuation ring, then $U=\Spec(K)$. If $k$ is perfect of characteristic $p$, then $C(k)$ is the ring $W(k)$ of $p$-typical Witt vectors with coefficients in $k$. We think of a $\BT$ over a scheme as an inductive system indexed by elements of $\mathbb N^*$ or alternatively as its limit $\textup{fppf}$ sheaf, see \cite{CCO}, Subsect. 3.3.2.

\section{Basic Extension Properties}\label{S2}

In this section we gather some basic extension properties of different nature that will be often used in what follows. We begin by recalling some standard properties of algebraic geometry.

\subsection{Reflexive sheaves}\label{SS21}

The proof of the following elementary fact is left to the reader.

\begin{fact}\label{F1} 
Let $\mathcal N$ be a normal noetherian scheme. Let $\mathcal U$ be an open dense subscheme of $\mathcal N$ which contains all codimension $1$ points (i.e., points whose local rings are discrete valuation rings). Let $\mathcal E$ and $\mathcal H$ be two finite flat group schemes over $\mathcal N$. Let $\mathcal V$ be a reflexive coherent $\mathcal O_{\mathcal N}$-module. Then the following four properties hold:

\medskip
{\bf (a)} If $\mathcal N$ is regular of dimension at most $2$, then the restriction from $\mathcal N$ to $\mathcal U$ induces an equivalence of categories from the category of coherent locally free $\mathcal O_{\mathcal N}$-modules to the category of coherent locally free $\mathcal O_{\mathcal U}$-modules.

\smallskip
{\bf (b)} If $\mathcal N$ is regular of dimension at most $2$, then the restriction from $\mathcal N$ to $\mathcal U$ induces an equivalence of categories from the category of finite flat group schemes over $\mathcal N$ to the category of finite flat group schemes over $\mathcal U$. 

\smallskip
{\bf (c)} The restriction homomorphisms $\textup{Hom}(\mathcal E,\mathcal H)\rightarrow\textup{Hom}({\mathcal U}\times_{\mathcal N} \mathcal E,{\mathcal U}\times_{\mathcal N} \mathcal H)$ and $H^0(\mathcal N,\mathcal V)\rightarrow H^0(\mathcal U,\mathcal V)$ are isomorphisms.

\smallskip
{\bf (d)} Let $f:\mathcal E\rightarrow \mathcal H$ be a homomorphism. Let $\mathcal W$ be the largest open subscheme of $\mathcal N$ such that the restriction $f_{\mathcal W}:\mathcal E_{\mathcal W}\rightarrow \mathcal H_{\mathcal W}$ of $f$ to a homomorphism over $\mathcal W$ is an isomorphism. Then $\mathcal W$ is the set of all points $z\in\mathcal N$ at which the fiber $f_z$ of $f$ is an isomorphism. Moreover, $\mathcal W$ contains all points of codimension at most $1$ in $\mathcal N$ if and only if we have $\mathcal W=\mathcal N$ (i.e., $f$ is an isomorphism if and only if its fibers over points of codimension at most $1$ in $\mathcal N$ are isomorphisms).
\end{fact}

For the second part of (c) above we refer to \cite{H}, Prop. 1.6.

\subsection{A refinement of a result of \cite{VZ3}}\label{SS22}

In this subsection we assume that $R$ is a discrete valuation ring of mixed characteristic $(0,p)$. Thus $K=\Frac(R)$ and $U=\Spec(K)$. Let $s\in\mathbb N$ be such that \cite{VZ3}, Thm. 1 holds for $R$; it has upper bounds in terms only on the absolute ramification index $e$ of $R$ (the best upper bound $s=[\log_p(\frac{pe}{p-1})]$ is obtained in  \cite{B}, Thm. A). If $H$ is a finite flat commutative group scheme of $p$ power order over $X$, for $n\in\mathbb N$ let $H\ceil{p^n}$ be the schematic closure of $H_K[p^n]$ in $H$; so $H\ceil{p^0}=H\ceil{1}$ is the trivial subgroup scheme of $H$. 

\medskip
The following lemma refines \cite{VZ3}, Cor. 4 which worked in a context in which $H^\prime$ is a $\BT_n$ over $X$.

\begin{lemma}\label{L1}
With $R$, $K$ and $s$ as above, let $n>2s$ be an integer. Let $H$ and $H^\prime$ be finite flat group scheme over $X$ such that their generic fibers $H_K$ and $H^\prime_K$ are $\BT_n$'s over $U=\Spec(K)$. For $i\in\{0,\ldots,n\}$ let $H_i:=H\ceil{p^i}$ and $H_i^\prime:=H^\prime\ceil{p^i}$. Then the following three properties hold:

\medskip
{\bf (a)} If $f:H\rightarrow H^\prime$ is a homomorphism such that $f_K:H_K\rightarrow H^\prime_K$ is an isomorphism, then the homomorphism $H_{n-s}/H_s\rightarrow H_{n-s}^\prime/H_s^\prime$ induced by $f$ is an isomorphism.

\smallskip
{\bf (b)} If $H_K\rightarrow H^\prime_K$ is a homomorphism (resp. is an isomorphism), then the homomorphism $H_K[p^{n-s}]/H_K[p^s]\rightarrow H_K^\prime[p^{n-s}]/H_K^\prime[p^s]$ induced by it extends to a homomorphism  (resp. to an isomorphism) $H_{n-s}/H_s\rightarrow H_{n-s}^\prime/H_s^\prime$.

\smallskip
{\bf (c)} If $n\ge 2s+2$, then $H_{s+t}/H_s$ is a $\BT_t$ over $X$ for all $t\in\{1,\ldots,n-2s\}$.
\end{lemma}

\noindent
{\bf Proof:} For $i\in\{0,\ldots,n-1\}$, let $\textup{gr}_{i+1}(H):=H_{i+1}/H_i$ and $\textup{gr}_{i+1}(H^\prime):=H_{i+1}^\prime/H_i^\prime$. To prove (a), we consider the commutative diagram
\[\xymatrixcolsep{8pc}\xymatrix{
H/H_{n-s-1} \ar[r]^{p^{n-s-1}} \ar[d] & H_{s+1} \ar[d] \\
H^\prime/H^\prime_{n-s-1} \ar[r]^{p^{n-s-1}} & H^\prime_{s+1} }\]
whose vertical arrows are induced by $f$. Let $\alpha:H/H_{n-s-1}\rightarrow H^\prime_{s+1}$ be the diagonal homomorphism of the diagram and let $\beta:H^\prime_{s+1}\rightarrow H/H_{n-s-1}$ be a homomorphism such that $\alpha\circ\beta$ and $\beta\circ\alpha$ are both the multiplication by $p^s$, cf. \cite{VZ3}, Thm. 1. The homomorphism $\beta$ induces a homomorphism $\textup{gr}_{s+1}(H^\prime)\rightarrow \textup{gr}_{n-s}(H)$ which is an inverse of the homomorphism 
$\textup{gr}_{n-s}(H)\rightarrow\textup{gr}_{s+1}(H^\prime)$ induced by $\alpha$. Factoring the isomorphism $\textup{gr}_{n-s}(H)\rightarrow \textup{gr}_{s+1}(H^\prime)$ we get that for $i\in\{s+1,\ldots,n-s\}$ the homomorphisms $\textup{gr}_{i+1}(H)\rightarrow \textup{gr}_{i+1}(H^\prime)$ induced by $f$ are isomorphisms. This implies that (a) holds. 

Using the schematic closure $\tilde H$ of the graph of $H_K\rightarrow H_K^\prime$ in the product $H\times_{\Spec(R)} H^\prime$, the two projections $\tilde H\rightarrow H$ and $\tilde H\rightarrow H^\prime$ induce homomorphisms $\tilde H_{n-s}/\tilde H_s\rightarrow H_{n-s}/H_s$ and $\tilde H_{n-s}/\tilde H_s\rightarrow H_{n-s}^\prime/H_s^\prime$, the first one being an isomorphism (resp. both of them being isomorphisms) (cf. (a)). From this (b) follows. 

To prove (c), from the proof of (a) applied to the identity automorphism of $H$, we get that $\textup{gr}_{n-s}(H)\rightarrow\textup{gr}_{s+1}(H)$ is an isomorphism, which implies that for $i\in\{s+1,\ldots,n-s-1\}$ all successive homomorphisms (which are isomorphisms over $U$) $\textup{gr}_{i+1}(H)\rightarrow \textup{gr}_i(H)$ are isomorphisms. Thus $H_{n-s}/H_s$ is an $\textup{fppf}$ sheaf flat over $\mathbb Z/p^{n-2s}\mathbb Z$. From \cite{Me}, Ch. 1, Def. 1.2 we get that (c) holds for $t=n-2s$. If $1\le t\le n-2s-1$, then $(H_{n-s}/H_s)[p^t]$ is a $\BT_t$ over $X$ and is equal to the flat closed subgroup scheme $H_{s+t}/H_s$ of $H_{n-s}/H_s$ as this is so over $U$.\endproof

\subsection{On Raynaud's complement to Tate's extension theorem}\label{SS23}

\medskip
We recall Raynaud's complement to Tate's extension theorem (see \cite{R}, Prop. 2.3.1) and include a new proof of it based on Lemma \ref{L1}.

\begin{prop}\label{P1} 
We assume $R$ is a discrete valuation ring of mixed characteristic $(0,p)$. Let $D_K$ be a $\BT$ over $K=\Frac(R)$. If for each $n\in\mathbb N^*$, $D_K[p^n]$ extends to a finite flat group scheme $E_n$ over $X$, then $D_K$ extends to a $\BT$ over $X$. 
\end{prop}

\noindent
{\bf Proof:} If $t\ge 2$ is an integer, then for $n:=2s+t\ge 2s+2$, the quotient group scheme $F_t:=E_n\ceil{p^{n-s}}/E_n\ceil{p^s}$ is a $\BT_t$ which extends $D_K[p^t]$ and does not depend on $E_n$ (see Lemma \ref{L1} (b) and (c)). Taking $E_n=E_{n+1}\ceil{p^n}$, we get that $F_t=F_{t+1}[p^t]$. Thus $F_1:=F_t[p]$ does not depend on $t$ and the inductive system $F_m$ indexed by $m\in\mathbb N^*$ is a $\BT$ over $X$ which extends $D_K$.\endproof

\medskip
In the proof of Theorem \ref{T1} (b) we will need the following application of Proposition \ref{P1}.

\begin{prop}\label{P2}
Let $R$ be a regular local ring of mixed characteristic $(0,p)$ and dimension $d\ge 2$. Let $D_U$ be a $\BT$ over $U$ (the punctured spectrum of $R$). Let $O$ be the discrete valuation ring which dominates $R$, which has the same field of fractions $K$ as $R$, and which is a local ring of the blow up of $X$ along its closed point (thus the residue field of $O$ is a field of rational functions in $d-1$ variables over $k$). Then the generic fiber $D_K$ of $D_U$ over $\Spec(K)$ extends uniquely (up to unique isomorphism) to a $\BT$ over $\Spec(O)$.
\end{prop}

\noindent
{\bf Proof:} 
For each $n\in\mathbb N^*$, from Lemma \ref{L2} below applied with $(H_K,H_U)=(D_K[p^n],D_U[p^n])$ we get that the $\BT_n$ $D_K[p^n]$ over $\Spec(K)$ extends to a finite flat group scheme over $\Spec(O)$. Based on this, the proposition follows from Proposition \ref{P1} applied over $O$.\endproof

\begin{lemma}\label{L2}
Let $R$, $d$, $U$, $O$, $K$, and $X$ be as in Proposition \ref{P2}. Let $H_K$ be a finite group scheme over $\Spec(K)$ which extends to a finite flat group scheme $H_U$ over $U$. Then $H_K$ extends to a finite flat group scheme over $\Spec(O)$. 
\end{lemma}

\noindent
{\bf Proof:} Let $R_1$ be a regular local ring of dimension $2$ which dominates $R$, which is dominated by $O$, and whose punctured spectrum $U_1$ is $X_1\times_X U$, where $X_1:=\Spec(R_1)$. For instance, we can take $R_1$ to be a local ring of the blow up of $X$ along a regular closed subscheme of it of dimension $1$. From Fact \ref{F1} (b) we get that $U_1\times_U H_U$ extends uniquely (up to unique isomorphism) to a finite flat group scheme $H_1$ over $X_1$. The pullback of $H_1$ to $\Spec(O)$ is a finite flat group scheme over $\Spec(O)$ that extends $H_K$.\endproof

\subsection{A descent lemma}\label{SS24}
In what follows we will often use the following general descent lemma.

\begin{lemma}\label{L3}
Let $R\rightarrow R^{\prime}$ be a faithfully flat extension of noetherian local rings with $\dim(R)=\dim(R^{\prime})$ and $\depth\;(R)\ge 2$. If $R^{\prime}$ is $p$-quasi-healthy, then $R$ is $p$-quasi-healthy.
\end{lemma}

\noindent
{\bf Proof:} Let $D_U$ be a $\BT$ over $U$. We have to show that $D_U$ extends (automatically uniquely up to unique isomorphism) to a $\BT$ $D$ over $X=\Spec(R)$. Let $X^{\prime}:=\Spec(R^{\prime})$ and let $f:X^{\prime}\rightarrow X$ be the morphism induced by the monomorphism $R\rightarrow R^{\prime}$. From \cite{Ma}, Thm. 15.1 we get that the closed fiber of $f$ has dimension zero, and therefore $U^{\prime}:=f^{-1}(U)$ is the punctured spectrum of $R^{\prime}$. Thus $D_{U^{\prime}}=U^{\prime}\times_U D_U$ extends to a $\BT$ $D_{X^{\prime}}$ over $X^{\prime}$.

For $n\in \mathbb N^*$, let $A_n:=H^0(U,D_U[p^n])$. We have to show that the following two properties hold:

\medskip\noindent
{\bf (i)} for each $n\in\mathbb N^*$, $A_n$ is a finitely generated projective $R$-module (so for $t\in\mathbb N$, the $t$-fold tensor power of $A_n$ over $R$ maps isomorphically to $\mathcal O(D_U[p^n]^t)$), and thus the group scheme structure of $D_U[p^n]$ endows $A_n$ with a structure of commutative and cocommutative Hopf $R$-algebra that defines a finite flat commutative group scheme $D_n$ over $X$ that extends $D_U[p^n]$);

\smallskip\noindent
{\bf (ii)} for all $n,m\in \mathbb N^*$ the complex $0\rightarrow D_m\rightarrow D_{n+m}\rightarrow D_n\rightarrow 0$ that extends the short exact sequence $0\rightarrow D_U[p^m]\rightarrow D_U[p^{n+m}]\rightarrow D_U[p^n]\rightarrow 0,$
is a short exact sequence of commutative finite flat group schemes over $X$.

\medskip
If these two properties hold, then we can take $D$ to be the inductive system $D_n$. To check that properties (i) and (ii) hold we can work
locally in the $\textup{fpqc}$ topology and therefore these two properties follow from the fact that for each $n\in \mathbb N^*$ we have $D_{X^{\prime}}[p^n]=\Spec(R{^{\prime}}\otimes_R A_n)= X^{\prime}\times_X D_n$.\endproof

\subsection{Extending vector bundles}\label{SS25}
The following criterion of extending vector bundles will be used repeatedly in connection to extending truncated $\BT$'s.

\begin{lemma}\label{L4}
We assume that $\depth(R)\ge 3$. Let $y\in R$ be such that $(y,p)$ is a regular sequence of $R$ (thus $R$ is of mixed characteristic $(0,p)$). Let $S$ be either $R$ or $R/(y^n)$ for some $n\in \mathbb N^*$. Let $\mathcal F$ be a vector bundle over the punctured spectrum $U_S$ of $S$. If $S=R$, let $\tau$ be $y$ and if $S=R/(y^n)$, let $\tau$ be $p$. For each $t\in \mathbb N^*$, let $\iota_t:U_S\cap \Spec(S/(\tau^t))\rightarrow \Spec(S/(\tau^t))$ be the restriction modulo $\tau^t$ of the open embedding $\iota:U_S\rightarrow\Spec(S)$. We assume that the following two conditions hold:

\medskip\noindent
{\bf (i)} for each $t\in \mathbb N^*$, the direct image $\iota_{t,*} \mathcal F_t$ of the restriction $\mathcal F_t$ of $\mathcal F$ to $U_S\cap \Spec(S/(\tau^t))$ is a vector bundle $\mathcal F_t^+$ over $\Spec(S/(\tau^t))$;

\smallskip\noindent
{\bf (ii)} if $S=R/(y^n)$, then $\depth(R)\ge 4$.

\medskip 
 Then the direct image $\iota_{*} \mathcal F$ is a vector bundle over $\Spec(S)$ whose restriction to $\Spec(S/(\tau^t))$ maps isomorphically to $\mathcal F_t^+$ for all $t\in\mathbb N^*$.
\end{lemma}

\noindent
{\bf Proof:} If $S=R$ is regular, then the lemma is a particular case of \cite{FC}, Ch. V, Sect. 6, Lem. 6.6 applied over $\widehat{S}$. If $S=R$ (resp. if $S=R/(y^n)$), then from the Auslander--Buchsbaum formula we get that $\depth(S)$ is $\depth(R)-1\ge 2$ (resp. $\depth(R)-2\ge 2$). Thus, regardless of what $S$ is, the lemma is a particular case of \cite{G2}, Exp. IX, Prop. 1.4 and Ex. 1.5 applied over $\widehat{S}$. \endproof

\subsection{Extending $\BT$'s}\label{SS26}

The next lemma is a natural application of Lemma \ref{L4} and is the essence of the inductive step of the induction we will use to prove Theorem \ref{T1} (b).

\begin{lemma}\label{L5}
Let $R$ be a local noetherian ring of mixed characteristic $(0,p)$ such that  $\depth(R)\ge 4$. We assume there exists $y\in R$ such that $(y,p)$ is a regular sequence of $R$. Let $S:=R/(y)$. Let $D_{U}$ be a $\BT$ over the punctured spectrum $U$ of $R$. If the restriction of $D_U$ to the punctured spectrum of $S$ extends to $\Spec(S)$, then $D_U$ extends to a $\BT$ over $X=\Spec(R)$. 
\end{lemma}

\noindent
{\bf Proof:} Based on the proof of Lemma \ref{L3} we can assume that $R$ is complete. For $t\in\mathbb N^*$ let 
$$X_t:=\Spec(R/(y)^t)$$ 
and $U_t:=U\cap X_t$. By induction on $t\in\mathbb N^*$ we will show that the restriction $D_{U_t}$ of $D_{U}$ to $U_t$ extends uniquely (up to unique isomorphism) to a $\BT$ $D_{X_t}$ over $X_t$. We already know that this holds for $t=1$. Let $\iota_{U_1}: D_{U_1}\rightarrow U_1\times_{X_1} D_{X_1}$ be the canonical isomorphism. 

The passage from $t$ to $t+1$ goes as follows. Assuming that $D_{X_t}$ exists, let $\iota_{U_t}: D_{U_t}\rightarrow U_t\times_{X_t} D_{X_t}$ be the canonical isomorphism. For $q\in \mathbb N^*$, let $D_{X_t,q}$ and $D_{U_t,q}$ be the reductions modulo $p^q$ of $D_{X_t}$ and $D_{U_t}$ (respectively); they are $\BT$'s over the reductions $X_{t,q}$ and $U_{t,q}$ of $X_t$ and $U_t$ (respectively) modulo $p^q$. From the Grothendieck--Messing deformation theory (see \cite{I}, Thm. 4.4 and Cor. 4.7 and \cite{Me}, Ch. V, Thm. 1.6) we get that the lifts of $D_{X_t,q}$ to $\BT$'s over $X_{t+1,q}$ are parametrized by the global sections of a torsor under the group of global sections of a coherent locally free $\mathcal O_{\Spec(S/p^qS)}$-module $\mathcal F_{t,q}$ (note that the ideal that defines the closed embedding $X_{t,q}\rightarrow X_{t+1,q}$ has a canonical trivial divided power structure and therefore this closed embedding has a canonical structure of a nilpotent divided power thickening). Similarly, the lifts of $D_{U_t,q}$ to $\BT$'s over $U_{t+1,q}$ are parametrized by a torsor under the group $H^0(U\cap \Spec(S/p^qS),\mathcal F_{t,q})$. As $(y,p)$ is a regular sequence of $R$, we have $\depth_S (S/p^qS)=\depth(R)-2\ge 2$ (cf. Auslander--Buchsbaum formula) and thus we can identify $H^0(U\cap \Spec(S/p^qS),\mathcal F_{t,q})=H^0(\Spec(S/p^qS),\mathcal F_{t,q})$. This implies that there exists a unique (up to unique isomorphism) $\BT$ $D_{X_{t+1},q}$ over $X_{t+1,q}$ that lifts compatibly both $D_{X_t,q}$ and $D_{U_{t+1},q}$. As $R/(y)^{t+1}$ is $p$-adically complete, there exists a unique (up to unique isomorphism) $\BT$ $D_{X_{t+1}}$ over $X_{t+1}$ which lifts the system $D_{X_{t+1},q}$, $q\in \mathbb N^*$ (cf. \cite{Me}, Ch. II, Lem. 4.16 which implies that the categories of $\BT$'s over $X_{t+1}$ and over the formal scheme which is the $p$-adic completion of $X_{t+1}$ are canonically equivalent). 

For each $n\in \mathbb N^*$, from Lemma \ref{L4} applied with 
$$\big(S,\tau,(\mathcal F_q,\mathcal F_q^+)_{q\in\mathbb N^*}\big)=\big(R/(y)^{t+1},p,(\mathcal O_{D_{U_{t+1},q}[p^n]},\mathcal O_{D_{X_{t+1},q}[p^n]})_{q\in\mathbb N^*}\big),$$ we get that $\mathcal O_{D_{U_{t+1}}[p^n]}$ is the restriction to $U_{t+1}$ of a vector bundle over $X_{t+1}$ whose restriction to $X_{t+1,q}$ is compatibly identified with $\mathcal O_{D_{X_{t+1},q}[p^n]}$ for all $q\in\mathbb N^*$, and thus (as $X_{t+1}$ is $p$-adically complete) this vector bundle can be functorially identified with $\mathcal O_{D_{X_{t+1}}[p^n]}$. We deduce an isomorphism $\iota_{U_{t+1}}[p^n]:D_{U_{t+1}}[p^n]\rightarrow U_{t+1}\times_{X_{t+1}} D_{X_{t+1}}[p^n]$ of $U_{t+1}$-schemes, and it must be a group scheme isomorphism. The $\iota_{U_{t+1}}[p^n]$'s glue together to define an isomorphism $\iota_{U_{t+1}}:D_{U_{t+1}}\rightarrow U_{t+1}\times_{X_{t+1}} D_{X_{t+1}}$. Therefore the restriction of $D_{X_{t+1}}$ to $U_{t+1}$ is canonically identified with $D_{U_{t+1}}$. This ends the induction.

From Lemma \ref{L4} applied with $(S,\tau)=(R,y)$ we similarly get that for each $n\in\mathbb N^*$ the locally free $\mathcal O_U$-module associated to $D_U[p^n]$ extends to a locally free $\mathcal O_X$-module whose reduction modulo each $(y)^t$ is the $\mathcal O_{X_t}$-module associated to $D_{X_t}[p^n]$ and which defines naturally a $\BT_n$ $D_n$ over $X$. The inductive system $D_n$ is a $\BT$ over $X$ which extends $D_U$ (and thus our notation matches); we note that it is also the unique (up to unique isomorphism) $\BT$ over $X$ which lifts compatibly $D_{X_t}$ for all $t\in \mathbb N^*$ (again cf. \cite{Me}, Ch. II, Lem. 4.16). Thus the lemma holds.\endproof 

\section{On the condition ($\natural$)}\label{S3}
In this section we assume that $R$ is regular of dimension $d\ge 1$ and mixed characteristic $(0,p)$ and study the condition ($\natural$) introduced in Theorem \ref{T1}. Let $C(k)$ be a Cohen ring with residue field $k$ which is a coefficient ring of $\widehat{R}$, as in the Cohen structure theorem for complete local noetherian rings of mixed characteristic (see \cite{G3}, Ch. 0, Subsect. 19.8 or \cite{Ma}, Subsect. 29). Thus $C(k)$ is a subring of $\widehat{R}$. Let $y_1,\ldots,y_d$ be a regular system of parameters of $\widehat{R}$. The natural $C(k)$-algebra homomorphism $\varrho:C(k)[[x_1,\ldots,x_d]]\rightarrow\widehat{R}$ that maps $x_i$ to $y_i$ is onto. This implies that we can identify 
$$\widehat{R}=C(k)[[x_1,\ldots,x_d]]/(h)$$ 
for some element $h\in\Ker(\varrho)\subset C(k)[[x_1,\ldots,x_d]]$ whose reduction $\bar h$ modulo $p$ is a non-zero element of the ideal $(x_1,\ldots,x_d)$ of $k[[x_1,\ldots,x_d]]$. The reduction $h_p$ of $h$ modulo the ideal $(x_1,\ldots,x_d)$ of $C(k)[[x_1,\ldots,x_d]]$ is an element of $C(k)$ with the property that $C(k)/(h_p)=\widehat{R}/(y_1,\ldots,y_d)=k$. Thus $(h_p)=(p)\subset C(k)$ and therefore $h_p$ is $p$ times a unit of $C(k)$. 

The ring $R/pR$ is regular if and only if $\bar h\notin (x_1,\ldots,x_d)^2$ and if and only if $\widehat{R}$ is isomorphic to $C(k)[[x_1,\ldots,x_{d-1}]]$. If $R/pR$ is regular, then there are $C(k)$-algebra epimorphisms $\theta:C(k)[[x_1,\ldots,x_d]]\rightarrow\widehat{R}$ under which the images of $x_1,\ldots,x_d$ in $\widehat{R}$ do not form a regular system of parameters of $\widehat{R}$; but each such $C(k)$-algebra epimorphism has a kernel generated by an element of $C(k)[[x_1,\ldots,x_d]]$ whose reduction modulo $p$ belongs to a regular system of parameters of $k[[x_1,\ldots,x_d]]$ and thus differs from $\bar h$ by a $k$-algebra automorphism of $k[[x_1,\ldots,x_d]]$. If $R/pR$ is not regular, then there exists no $C(k)$-algebra epimorphism $\theta:C(k)[[x_1,\ldots,x_d]]\rightarrow\widehat{R}$ under which the images of $x_1,\ldots,x_d$ in $\widehat{R}$ do not form a regular system of parameters of $\widehat{R}$.

If $C^{\prime}(k)$ is another Cohen ring which is a coefficient ring of $\widehat{R}$ and if $\varrho^{\prime}:C^{\prime}(k)[[x_1,\ldots,x_d]]\rightarrow\widehat{R}$ is a $C^{\prime}(k)$-algebra epimorphism such that the images of $x_1,\ldots,x_d$ form a regular system of parameters of $\widehat{R}$, then there exists an isomorphism $\nu:C(k)[[x_1,\ldots,x_d]]\rightarrow C^{\prime}(k)[[x_1,\ldots,x_d]]$ such that we have $\varrho^{\prime}=\varrho\circ\nu$. To check this, up to a $C^{\prime}(k)$-algebra automorphism of $C^{\prime}(k)[[x_1,\ldots,x_d]]$ preserving the ideal $(x_1,\ldots,x_d)$, we can assume that $\varrho^{\prime}(x_i)=y_i$. As the homomorphism $\mathbb Z_p\rightarrow C(k)$ is formally smooth, there exists a homomorphism $\nu_p:C(k)\rightarrow C^{\prime}(k)[[x_1,\ldots,x_d]]$ whose composite with $\varrho^{\prime}$ is the inclusion $C(k)\rightarrow\widehat{R}$. Thus we can take $\nu$ such that it extends $\nu_p$ and maps $x_i$ into $x_i$ for all $i\in\{1,\ldots,d\}$.

From the last two paragraphs we get that the ideal $(\bar h)$ of $k[[x_1,\ldots,x_d]]$ coming from a presentation of $\widehat{R}$ as $C(k)[[x_1,\ldots,x_d]]/(h)$ is uniquely determined by $R$ up to automorphisms of $k[[x_1,\ldots,x_d]]$ inducing the identity on the residue field $k$, and $\bar h$ is uniquely determined up to units and such automorphisms. We say that $\bar h$ is a Cohen element of $R$. Thus condition ($\natural$) for $R$ means that one Cohen element of $R$ does not belong to the ideal $(x_1^p,\ldots,x_d^p)+(x_1,\ldots,x_d)^{2p-2}$ of $k[[x_1,\ldots,x_d]]$, equivalently all Cohen elements of $R$ have this property.

Recall (\cite{Ma}, Exerc. 14.5) that the multiplicity of $R/pR$ is the order of vanishing of $p$, i.e., it is the positive integer $e$ such that $p\in \mathfrak m^e\setminus \mathfrak m^{e+1}$, where $\mathfrak m$ is the maximal ideal of $R$. As $\widehat{R}/p\widehat{R}\simeq k[[x_1,\ldots,x_d]]/(\bar h)$, $e$ is also the order of vanishing of $\bar h$. In particular, we have: 

\medskip\noindent
{\bf (i)} if $e<p$, then condition ($\natural$) holds;

\smallskip\noindent
{\bf (ii)} if condition ($\natural$) holds, then $e\le 2p-3$ for $d\ge 2$ and $e\le p-1$ for $d=1$.

\begin{fact}\label{F4}
Suppose $\widehat{R}$ is presented as $C(k)[[x_1,\ldots,x_d]]/(h)$. Let $l$ be a field extension of $k$ and let $C(k)\rightarrow C(l)$ be a monomorphism between Cohen rings which modulo $p$ is the inclusion $k\rightarrow l$, cf. \cite{C}, Cor. 1. Then the $R$-algebra $R^{\prime}:=C(l)[[x_1,\ldots,x_d]]/(h)$ is faithfully flat, $R^{\prime}$ is regular, and condition ($\natural$) holds for $R$ if and only if it holds for $R^{\prime}$.
\end{fact}

\noindent
{\bf Proof:} As the $C(k)$-algebra $C(l)$ is faithfully flat, the $C(k)[x_1,\ldots,x_d]$-algebra $C(l)[x_1,\ldots,x_d]$ is also faithfully flat. Taking completions of localizations with respect to maximal ideals $(p,x_1,\ldots,x_d)$ we get that the $C(k)[[x_1,\ldots,x_d]]$-algebra $C(l)[[x_1,\ldots,x_d]]$ is flat, see \cite{Ma}, Thm. 22.4 (i). Hence $\widehat{R}\rightarrow R^{\prime}$ is flat; being local, it is faithfully flat. The same holds for $R\rightarrow R^{\prime}$, and $R^{\prime}$ is regular by \cite{Ma}, Thm. 23.7. The last assertion is clear.\endproof

\subsection{Taking slices}\label{SS31}

In this subsection we study for $d\ge 3$ the behavior of the condition ($\natural)$ under restrictions to quotient rings of $\widehat{R}$ which are regular of dimension $d-1$. We have the following general and abstract lemma on homogeneous polynomials which also holds for $d=2$.

\begin{lemma}\label{L6} Let $l$ be a field of characteristic $p$. Let $d\ge 2$ be an integer and let $c\in\{1,\ldots,(d-1)(p-1)\}$. Let $\mathbb I_c$ be the subset of $\{0,\ldots,p-1\}^d$ formed by all those $d$-tuples whose entries sum up to $c$. We consider a non-zero homogeneous polynomial 
$$\bar h_c(x_1,\ldots,x_d)=\sum_{(i_1,\ldots,i_d)\in\mathbb I_c} \delta_{i_1,\ldots,i_d}\prod_{t=1}^d x_t^{i_t}\in l[x_1,\ldots,x_d]$$ 
of degree $c$. Then there exists an open dense subscheme $\mathbb O$ of the affine scheme $\mathbb A^{d-1}_l$ such that for each $d-1$-tuples $(v_1,\ldots,v_{d-1})\in \mathbb O(l)\subset \mathbb A^{d-1}_l(l)=l^{d-1}$ the following homogeneous polynomial 
$$\bar q_c(x_1,\ldots,x_{d-1})=\bar h_c(x_1,\ldots,x_{d-1},\sum_{t=1}^{d-1} v_tx_t)$$ has a non-zero image in $l[x_1,\ldots,x_{d-1}]/(x_1^p,\ldots,x_{d-1}^p)$.
\end{lemma}

\noindent
{\bf Proof:} For $q\in\mathbb N\cap \{c-p+1,\ldots,c\}$ let $\mathbb J_q^{d-1}$ be the subset of $\{0,\ldots,p-1\}^{d-1}$ formed by all those $d-1$-tuples whose entries sum up to $q$. For $(j_1,\ldots,j_{d-1})\in\mathbb J_c^{d-1}$, the coefficient of $\prod_{t=1}^{d-1} x_t^{j_t}$ in $\bar q_c$ is a sum of the form 
$$\sum_{q\in\mathbb N\cap \{c-p+1,\ldots,c\}} \sum_{(i_1,\ldots,i_{d-1})\in\mathbb J_q^{d-1},i_1\le j_1,\ldots,i_{d-1}\le j_{d-1}} \gamma^{j_1,\ldots,j_{d-1}}_{i_1,\ldots,i_{d-1}}\delta_{i_1,\ldots,i_{d-1},c-q}\prod_{t=1}^{d-1} v_t^{j_t-i_t},$$
where each $\gamma^{j_1,\ldots,j_{d-1}}_{i_1,\ldots,i_{d-1}}$ is a multinomial coefficient that divides $(c-q)!$ and thus is a non-zero element of $l$. 

If $\bar h_c$ is a polynomial only in the variables $x_1,\ldots,x_{d-1}$, then $\bar q_c=\bar h_c$ does not depend on the $d-1$-tuple $(v_1,\ldots,v_{d-1})$ and has a non-zero image in $l[x_1,\ldots,x_{d-1}]/(x_1^p,\ldots,x_{d-1}^p)$; thus in this case we can take $\mathbb O:=\mathbb A_l^{d-1}$.

Therefore we can assume that $\bar h_c$ is not a polynomial only in the variables $x_1,\ldots,x_{d-1}$, i.e., there exists a non-zero coefficient $\delta_{i_1,\ldots,i_d}$ of $\bar h_c$ with $i_d\ge 1$ or equivalently with $\sum_{t=1}^{d-1} i_t<c$. As $c\le (d-1)(p-1)$, if for each $t\in\{1,\ldots,d-1\}$ we consider an element $j_t\in\{i_t,\ldots,p-1\}$ such that we have $\sum_{t=1}^{d-1} j_t=c$, then the coefficient of $\prod_{t=1}^{d-1} x_t^{j_t}$ in $\bar q_c$ is a non-zero polynomial $\bar q_{c,j_1,\ldots,j_{d-1}}(v_1,\ldots,v_{d-1})\in l[v_1,\ldots,v_{d-1}]$. We consider the open dense subscheme $\mathbb O$ of $\mathbb A_l^{d-1}$ such that for a $d-1$-tuple $(v_1,\ldots,v_{d-1})\in\mathbb A^{d-1}_l(l)=l^{d-1}$ we have $\bar q_{c,j_1,\ldots,j_{d-1}}(v_1,\ldots,v_{d-1})\neq 0$ if and only if $(v_1,\ldots,v_{d-1})\in\mathbb O(l)$. Thus if $(v_1,\ldots,v_{d-1})\in\mathbb O(l)$, the coefficient of $\prod_{t=1}^{d-1} x_t^{j_t}$ in $\bar q_c$ is non-zero and therefore the image of $\bar q_c$ in $l[x_1,\ldots,x_{d-1}]/(x_1^p,\ldots,x_{d-1}^p)$ is non-zero.\endproof 

\begin{prop}\label{P3} 
We assume that $d\ge 3$, that $k$ is an infinite field, that $R=\widehat{R}$ is complete, and that condition ($\natural$) holds for $R$. Then for each $d_1\in\{2,\ldots,d-1\}$ there exists a closed regular subscheme $\Spec(S)$ of $X=\Spec(R)$ of dimension $d_1$ and mixed characteristic $(0,p)$ such that condition ($\natural$) holds for $S$.
\end{prop}

\noindent
{\bf Proof:} Proceeding by induction on $d\ge 3$, it suffices to consider the case when $d_1=d-1$. We write 
$\bar h=\bar h_0+\sum_{i=1}^{2p-3} \bar h_i,$
where $\bar h_0$ belongs to the ideal $(x_1^p,\ldots,x_d^p)+(x_1,\ldots,x_d)^{2p-2}$ of $k[[x_1,\ldots,x_d]]$ and where each $\bar h_i$ is either $0$ or a homogeneous polynomial in the variables $x_1,\ldots,x_d$ of degree $i$. Let $c\in\{1,\ldots,2p-3\}$ be the smallest integer such that $\bar h_c$ is non-zero (it exists as condition ($\natural$) holds for $R$). The field $k$ is infinite and we have $c\le 2p-3< (d-1)(p-1)$ as $d\ge 3$. Thus, from Lemma \ref{L6} applied with $l=k$ we get that there exists a $d-1$-tuple $(v_1,\ldots,v_{d-1})\in\mathbb O(k)\subset k^{d-1}$ such that the homogeneous polynomial $\bar h_c(x_1,\ldots,x_{d-1},\sum_{t=1}^{d-1} v_tx_t)$ of degree $c$ has a non-zero-image in $k[[x_1,\ldots,x_{d-1}]]/(x_1^p,\ldots,x_{d-1}^p)$. This implies that $\bar h(x_1,\ldots,x_{d-1},\sum_{t=1}^{d-1} v_tx_t)$ does not belong to the ideal $(x_1^p,\ldots,x_{d-1}^p)+(x_1,\ldots,x_{d-1})^{2p-2}$ of $k[[x_1,\ldots,x_{d-1}]]$ and therefore condition ($\natural$) holds for the regular local ring 
$$S:=R/(x_d-\sum_{t=1}^{d-1} w_ix_i)=C(k)[[x_1,\ldots,x_{d-1}]]/(h(x_1,\ldots,x_{d-1},\sum_{t=1}^{d-1} w_tx_t))$$ 
of dimension $d-1$ and mixed characteristic $(0,p)$, where $(w_1,\ldots,w_{d-1})\in C(k)^{d-1}$ is such that its reduction modulo $p$ is $(v_1,\ldots,v_{d-1})\in k^{d-1}$.\endproof

\subsection{Functorial properties}\label{SS32}

In this subsection we study the behavior of the condition ($\natural$) under local homomorphisms. Let $\varphi:R\rightarrow R^{\prime}$ be a local homomorphism between regular local rings of mixed characteristic $(0,p)$. Let $d^{\prime}\ge 1$ be the dimension of $R^{\prime}$. Let $k^{\prime}$ be the residue field of $R^{\prime}$ and let $C(k^{\prime})$ be a Cohen ring which is a coefficient ring of $\widehat{R^{\prime}}$. Let $y_1^{\prime},\ldots,y_{d^{\prime}}^{\prime}$ be a regular system of parameters of $\widehat{R^{\prime}}$. The natural $C(k^{\prime})$-algebra homomorphism $\varrho^{\prime}:C(k^{\prime})[[x_1^{\prime},\ldots,x_{d^{\prime}}^{\prime}]]\rightarrow\widehat{R^{\prime}}$ that maps $x_i^{\prime}$ to $y_i^{\prime}$ is onto and we identify 
$$\widehat{R^{\prime}}=C(k^{\prime})[[x_1^{\prime},\ldots,x_{d^{\prime}}^{\prime}]]/(h^{\prime})$$ 
for some element $h^{\prime}\in C(k^{\prime})[[x_1^{\prime},\ldots,x_{d^{\prime}}^{\prime}]]$. 
As the homomorphism $\mathbb Z_p\rightarrow C(k)$ is formally smooth, the composite homomorphism $C(k)\rightarrow \widehat{R}\rightarrow \widehat{R^{\prime}}$ admits a lift to a homomorphism 
$\Phi_p: C(k)\rightarrow C(k^{\prime})[[x_1^{\prime},\ldots,x_{d^{\prime}}^{\prime}]]$. We consider a homomorphism $\Phi:C(k)[[x_1,\ldots,x_d]]\rightarrow C(k^{\prime})[[x_1^{\prime},\ldots,x_{d^{\prime}}]]$ that extends $\Phi_p$, that maps each $x_i$ into the ideal $(x_1^{\prime},\ldots, x_{d^{\prime}}^{\prime})$ of $C(k^{\prime})[[x_1^{\prime},\ldots,x_{d^{\prime}}]]$, and such that we have a commutative diagram:
\[\xymatrix{
C(k)[[x_1,\ldots,x_d]] \ar[r]^{\Phi} \ar[d]^{\varrho} & C(k^{\prime})[[x_1^{\prime},\ldots,x_{d^{\prime}}]] \ar[d]^{\varrho^{\prime}} \\
\widehat{R} \ar[r]^{\widehat{\varphi}} & \widehat{R^{\prime}}.
}\]

\begin{prop}\label{P4} 
With $\varphi$ and $\Phi$ as above, the following three properties hold:

\medskip
{\bf (a)} We have an equality $(\Phi(h))=(h^{\prime})$ of ideals of $C(k^{\prime})[[x_1^{\prime},\ldots,x_{d^{\prime}}^{\prime}]]$.

\smallskip
{\bf (b)} If the condition $(\natural)$ holds for $R^{\prime}$, then it also holds for $R$.

\smallskip
{\bf (c)} Let $X^{\prime}:=\Spec(R^{\prime})$. We assume that the morphism $X^{\prime}\rightarrow X$ defined by $\varphi$ is flat and its fiber over the closed point of $X$ is regular (i.e., the local ring $R^{\prime}/(\varphi(y_1),\ldots,\varphi(y_d))R^{\prime}$ is regular). Then the condition $(\natural)$ holds for $R$ if and only if it holds for $R^{\prime}$.
\end{prop}

\noindent
{\bf Proof:} As $\Phi$ induces $\widehat{\varphi}:\widehat{R}\rightarrow\widehat{R^{\prime}}$, we have $\Phi(h)\in (h^{\prime})$. We consider the unique element $u\in C(k^{\prime})[[x_1^{\prime},\ldots,x_{d^{\prime}}^{\prime}]]$ such that we have $\Phi(h)=h^{\prime}u$. To prove (a) it suffices to show that $u$ is a unit of $C(k^{\prime})[[x_1^{\prime},\ldots,x_{d^{\prime}}^{\prime}]]$. If $h^{\prime}_p$ and $u_p$ are the reductions modulo the ideal $(x_1^{\prime},\ldots,x_{d^{\prime}}^{\prime})$ of $C(k^{\prime})[[x_1^{\prime},\ldots,x_{d^{\prime}}^{\prime}]]$, then the equality $\Phi(h)=h^{\prime}u$ modulo the ideal $(x_1^{\prime},\ldots,x_{d^{\prime}}^{\prime})$ becomes an identity $\varphi_p(h_p)=h^{\prime}_pu_p$ between elements of $C(k^{\prime})$, where $\varphi_p:C(k)\rightarrow C(k^{\prime})=C(k^{\prime})[[x_1^{\prime},\ldots,x_{d^{\prime}}^{\prime}]]/(x_1^{\prime},\ldots,x_{d^{\prime}}^{\prime})$ is induced by $\Phi_p$. We know that $h_p$ and $h^{\prime}_p$ are $p$ times units of $C(k)$ and $C(k^{\prime})$ (respectively), see the beginning of Section \ref{S3}. From the last two sentences we get that $u_p$ is a unit of $C(k^{\prime})$ and therefore $u$ is a unit of $C(k^{\prime})[[x_1^{\prime},\ldots,x_{d^{\prime}}^{\prime}]]$. Thus (a) holds.

To prove (b) let $\bar h^{\prime}$ and $\bar\Phi: k[[x_1,\ldots,x_d]\rightarrow k^{\prime}[[x_1^{\prime},\ldots,x_{d^{\prime}}^{\prime}]]$ be the reductions modulo $p$ of $h^{\prime}$ and $\Phi$ (respectively). The elements $\bar\Phi(\bar h)$ and $\bar h^{\prime}$ differ by a unit, cf. (a). Thus if $\bar h$ belongs to the ideal $(x_1^p,\ldots,x_d^p)+(x_1,\ldots,x_d)^{2p-2}$ of $k[[x_1,\ldots,x_d]]$, then $\bar\Phi(\bar h)$ and $\bar h^{\prime}$ belong to the ideal of $k^{\prime}[[x_1^{\prime},\ldots,x_{d^{\prime}}^{\prime}]]$ generated by $\bar\Phi((x_1^p,\ldots,x_d^p))+\bar\Phi((x_1,\ldots,x_d)^{2p-2})$ and thus to the ideal $((x_1^{\prime})^p,\ldots,(x_{d^{\prime}}^{\prime})^p)+(x_1^{\prime},\ldots,x_{d^{\prime}}^{\prime})^{2p-2}$ of $k^{\prime}[[x_1^{\prime},\ldots,x_d^{\prime}]]$ and therefore the condition $(\natural)$ does not holds for $R^{\prime}$. From this (b) follows. 

Based on (b), to prove (c) it suffices to show that if the condition $(\natural)$ holds for $R$, then it also holds for $R^{\prime}$. To ease notation we can assume that $R$ and $R^{\prime}$ are complete; thus $\varphi=\widehat{\varphi}$, $y_1,\ldots,y_d\in R$, and $y^{\prime}_1,\ldots,y_{d^{\prime}}^{\prime}\in R^{\prime}$. Let $c$ be the dimension of the regular local ring $\bar R^{\prime}=R^{\prime}/(\varphi(y_1),\ldots,\varphi(y_d))R^{\prime}$ of characteristic $p$. As $\varphi$ is flat, we have $d^{\prime}=d+c$ (cf. \cite{Ma}, Thm. 15.1). We can assume that $y_{d+1}^{\prime},\ldots,y_{d+c}^{\prime}$ map into a regular system of parameters of $\bar R^{\prime}$. Let $R^{\prime}_1=R^{\prime}/(y_{d+1}^{\prime},\ldots,y_{d+c}^{\prime})$; it is a regular local ring of dimension $d=d^{\prime}-c$. The composite homomorphism $\varphi_1:R\rightarrow R^{\prime}_1$ between local rings of dimension $d$ is such that $R_1^{\prime}/(\varphi(y_1),\ldots,\varphi(y_d))=\Spec(k^{\prime})$. This implies that $\varphi_1$ is flat, cf. \cite{Ma}, Thm. 23.1. Thus $R_1^{\prime}$ has mixed characteristic $(0,p)$. If condition ($\natural$) holds for $R^{\prime}_1$, then it also holds for $R^{\prime}$ (cf. (b)) and therefore by replacing $R^{\prime}$ with $R^{\prime}_1$ we can assume that $d^{\prime}=d$, that $c=0$, and therefore that we have $y_i^{\prime}=\varphi(y_i)$ for all $i\in\{1,\ldots,d\}$. We can also assume that we have $x_i^{\prime}=\Phi(x_i)$ for all $i\in\{1,\ldots,d\}$, that $\bar\Phi:k[[x_1,\ldots,x_d]\rightarrow k^{\prime}[[x_1^{\prime},\ldots,x_{d^{\prime}}^{\prime}]]$ can be identified with the canonical inclusion $k[[x_1,\ldots,x_d]\rightarrow k^{\prime}[[x_1,\ldots,x_d]]$, and that under this identification we have $\bar h=\bar h^{\prime}$. As $\bar h$ does not belong to the ideal $(x_1^p,\ldots,x_d^p)+(x_1,\ldots,x_d)^{2p-2}$ of $k[[x_1,\ldots,x_d]]$, it does not belong to the ideal $(x_1^p,\ldots,x_d^p)+(x_1,\ldots,x_d)^{2p-2}$ of $k^{\prime}[[x_1,\ldots,x_d]]$ and thus the condition $(\natural)$ holds for $R^{\prime}$. Thus (c) holds.\endproof

\subsection{A generization property}\label{SS33}

\begin{prop}\label{P5} 
We assume that $R$ is a regular local ring of dimension $d\ge 2$ and mixed characteristic $(0,p)$ for which condition ($\natural$) holds. Let $R_1$ be a local ring of $X=\Spec(R)$ which is of mixed characteristic $(0,p)$ and dimension $d_1\in\{1,\ldots,d-1\}$. Then condition ($\natural$) holds for $R_1$.
\end{prop}

\noindent
{\bf Proof:} Let $k_1$ be the residue field of $R_1$. Let $\mathfrak p_1$ be the prime ideal of $R$ such that $R_1=R_{\mathfrak p_1}$. We consider two cases as follows.

\medskip
{\bf Case 1.} We first consider the case when $R=\widehat{R}$ is complete. We fix an identification $R=C(k)[[x_1,\ldots,x_d]]/(h)$, where $h\in C(k)[[x_1,\ldots,x_d]]$ is such that its reduction $\bar h$ modulo 
$p$ does not belong to the ideal $(x_1^p,\ldots,x_d^p)+(x_1,\ldots,x_d)^{2p-2}$ of $k[[x_1,\ldots,x_d]]$. Let $\mathfrak p:=\varrho^{-1}(\mathfrak p_1)\in\Spec(C(k)[[x_1,\ldots,x_d]])$. Let $\mathfrak n:=\mathfrak p C(k)[[x_1,\ldots,x_d]]_{\mathfrak p}$ be the maximal ideal of $C(k)[[x_1,\ldots,x_d]]_{\mathfrak p}$ and let $\mathfrak n^{(p)}$ be the $\mathfrak n$-primary ideal of $C(k)[[x_1,\ldots,x_d]]_{\mathfrak p}$ generated by $p$ and by $p$-th powers of elements of $\mathfrak n$. Let $\mathfrak A:=C(k)[[x_1,\ldots,x_d]]_{\mathfrak p}/\mathfrak n^{(p)}$; it is a local artinian ring which contains $k$, which has $\mathfrak j:=\mathfrak n/\mathfrak n^{(p)}$ as its maximal ideal, and which as a $k$-algebra is isomorphic to $k_1[y_1,\ldots,y_{d_1}]/(y_1^p,\ldots,y_{d_1}^p)$. 

As $\bar h$ does not belong to the ideal $(x_1^p,\ldots,x_d^p)+(x_1,\ldots,x_d)^{2p-2}$ of the ring $k[[x_1,\ldots,x_d]]$, we can speak about the smallest integer $n\in\{1,\ldots,2p-3\}$ such that $\bar h$ has in its writing a non-zero term of the form $v\prod_{j=1}^d x_j^{i_j}$ with each $i_j\in\{0,\ldots,p-1\}$ such that $n=\sum_{j=1}^d i_j$ and with $v\in k\setminus\{0\}$. This implies that there exist $n$ derivations $\partial_1,\ldots,\partial_n$ of $C(k)[[x_1,\ldots,x_d]]$ which are $C(k)$-linear, which belong to the set $\{\frac{\partial}{\partial x_1},\ldots,\frac{\partial}{\partial x_d}\}$, and for which the 
element $(\partial_1\circ\cdots \circ \partial_n)(h)$ is a unit of $C(k)[[x_1,\ldots,x_d]]$ (e.g., we can take $\partial_1=\cdots =\partial_{i_1}:=\frac{\partial}{\partial x_1}$, $\partial_{i_1+1}=\cdots =\partial_{i_1+i_2}:=\frac{\partial}{\partial x_2},\ldots,\partial_{i_1+\ldots+i_{d-1}+1}=\cdots=\partial_n=:\frac{\partial}{\partial x_d}$). For $i\in\{1,\ldots,n\}$ the derivation of $C(k)[[x_1,\ldots,x_d]]_\mathfrak p$ that extends $\partial_i$ induces a derivation $\bar\partial_i:\mathfrak A\rightarrow \mathfrak A$. 

Let $\bar h_{(p)}$ be the image of $\bar h$ in $\mathfrak A$. As $(\partial_1\circ\cdots \circ \partial_n)(h)$ is a unit of $C(k)[[x_1,\ldots,x_d]]$, the element $(\bar\partial_1\circ\cdots \circ \bar\partial_n)(\bar h_{(p)})$ is a unit of $\mathfrak A$. We show that the assumption that $\bar h_{(p)}\in \mathfrak j^{2p-2}$ leads to a contradiction. By a decreasing induction on $i\in\{1,\ldots,n\}$ we get that $(\bar\partial_i\circ\cdots \circ \bar\partial_n)(\bar h_{(p)})\in \mathfrak j^{2p-2-n+i-1}$. Thus the element $(\bar\partial_1\circ\cdots \circ \bar\partial_n)(\bar h_{(p)})$ belongs to $\mathfrak j^{2p-2-n}$. From this and the inequality $2p-2-n\ge 1$ we get that $(\bar\partial_1\circ\cdots \circ \bar\partial_n)(\bar h_{(p)})\in \mathfrak j$ is not a unit of $\mathfrak A$. Contradiction. Thus $\bar h_{(p)}\notin\mathfrak j^{2p-2}$. This implies that $h$ does not belong to the ideal $\mathfrak n^{(p)}+\mathfrak n^{2p-2}$ of $C(k)[[x_1,\ldots,x_d]]_{\mathfrak p}$.

Let $C(k_1)$ be a Cohen ring which is a coefficient ring of the completion $\reallywidehat{C(k)[[x_1,\ldots,x_d]]_{\mathfrak p}}$. We fix a $C(k_1)$-algebra isomorphism 
$$\reallywidehat{C(k)[[x_1,\ldots,x_d]]_{\mathfrak p}}\simeq C(k_1)[[y_1,\ldots,y_{d_1}]]$$ 
which lifts the $k_1$-algebra isomorphism $\mathfrak A/\mathfrak j\simeq k_1[y_1,\ldots,y_{d_1}]/(y_1^p,\ldots,y_{d_1}^p)$. The completion $\widehat{R_1}$ of $R_1$ is canonically identified with $\reallywidehat{C(k)[[x_1,\ldots,x_d]]_{\mathfrak p}}/(h)$ and thus with $C(k_1)[[y_1,\ldots,y_{d_1}]]/(h)$. From the last sentence of the previous paragraph we get that the image of $\bar h$ in $(y_1,\ldots,y_{d_1})$ does not belong to the ideal $(y_1^p,\ldots,y_{d_1}^p)+(y_1,\ldots,y_{d_1})^{2p-2}$ of $k_1[[y_1,\ldots,y_{d_1}]]$. Thus condition ($\natural$) holds for $R_1$ if $R=\widehat{R}$ is complete.

\medskip
{\bf Case 2.} We consider the case when $R$ is not complete. Let $\mathfrak q$ be a prime ideal of $\widehat{R}$ which is minimal over $\mathfrak p_1\widehat{R}$. We have a canonical homomorphism $R_1=R_{\mathfrak p_1}\rightarrow\widehat{R}_{\mathfrak q}$ of local rings. From Case 1 applied to $\widehat{R}$ we get that condition ($\natural$) holds for $\widehat{R}_{\mathfrak q}$. Thus condition ($\natural$) also holds for $R_1$, by Proposition \ref{P4} (b).\endproof

\section{On closed subschemes of projective spaces}\label{S4}
Let $l$ be a field. Let $N\in\mathbb N^*$. Let $\mathcal C$ be a closed connected subscheme (thus nonempty) of the projective space $\mathbb P^N_l$. Let $\reallywidehat{\mathbb P^N_l}$ be the formal completion of $\mathbb P^N_l$ along $\mathcal C$. Let $\kappa$ be the field of rational functions of $\mathbb P^N_l$ and let $\kappa^{\prime}$ be the ring of formal-rational functions on $\reallywidehat{\mathbb P^N_l}$. We recall from \cite{HM}, Def. (2.9.3) that $\mathcal C$ is said to be G3 in $\mathbb P^N_l$ if the homomorphism $\kappa\rightarrow\kappa^{\prime}$ associated to the natural morphism $\chi:\reallywidehat{\mathbb P^N_l}\rightarrow \mathbb P^N_l$ of locally ringed spaces is an isomorphism. For instance, from the proof of \cite{HM}, Thm. (3.3) we get that if $l$ is an algebraically closed field and if $\mathcal C$ is irreducible of positive dimension, then $\mathcal C$ is (universally) G3 in $\mathbb P^N_l$. In Section \ref{S5} we will need the following general result of formal algebraic geometry.

\begin{lemma}\label{L7}
We assume that $\mathcal C$ is a geometrically connected scheme of positive dimension which is $G3$ in $\mathbb P^N_l$. Let $\mathcal U$ be an open subscheme of $\mathbb P^N_l$ which contains $\mathcal C$. Let $\mathcal V$ be a torsion free coherent sheaf on $\mathcal U$. Let $\reallywidehat{\mathcal V}:=\chi^*(\mathcal V)$ be the coherent sheaf on $\reallywidehat{\mathbb P^N_l}$ which is the natural pullback of $\mathcal V$. Then there exists an open subscheme $\mathcal W$ of $\mathcal U$ which contains $\mathcal C$ and such that the natural pullback homomorphism 
$$\rho_{\mathcal W}:H^0(\mathcal W,\mathcal V)\rightarrow H^0(\reallywidehat{\mathbb P^N_l},\reallywidehat{\mathcal V})$$ 
between global sections is an isomorphism between finite dimensional $l$-vector spaces. Moreover, if $\mathcal V$ is a reflexive $\mathcal O_{\mathcal U}$-module, then we can take $\mathcal W=\mathcal U$. 
\end{lemma}

\noindent
{\bf Proof:} Using a meromorphic basis of $\mathcal V$ and the assumption that $\kappa=\kappa^{\prime}$, one sees that the finite dimensional $\kappa$-vector space $V$ of meromorphic sections of $\mathcal V$ maps isomorphically to the space of meromorphic sections of $\reallywidehat{\mathcal V}$. Thus, as $\mathcal V$ is torsion free, each homomorphism $\rho_{\mathcal W}$ is injective and we have natural inclusions $H^0(\mathcal U,\mathcal V)\subset H^0(\reallywidehat{\mathbb P^N_l},\reallywidehat{\mathcal V})\subset V$. 

We will first show that for each $\xi\in H^0(\reallywidehat{\mathbb P^N_l},\reallywidehat{\mathcal V})\subset V$ there exists an open subscheme $\mathcal W_{\xi}$ of $\mathcal U$ which contains $\mathcal C$ and such that we have $\xi\in H^0(\mathcal W_{\xi},\mathcal V)$. 

Let $\mathcal O$ be a local ring of $\mathbb P^N_l$ at a closed point $z$ of $\mathcal C$; it is a subring of $\kappa$. Let $\reallywidehat{\mathcal O}$ be the completion of $\mathcal O$ and let $\mathcal M$ be the $\mathcal O$-submodule of $V$ defined by the global sections of $\mathcal V$ over $\Spec(\mathcal O)$; we can identify $V=\kappa\otimes_{\mathcal O} \mathcal M$. 

As the homomorphism $\mathcal O\rightarrow \reallywidehat{\mathcal O}$ is faithfully flat, it is well-known that we have $\kappa\cap\reallywidehat{\mathcal O}=\mathcal O$ (the intersection being taken inside the field of fractions of $\reallywidehat{\mathcal O}$). We recall the argument for this. For $a\in \kappa\cap\reallywidehat{\mathcal O}$, we have inclusions of $\mathcal O$-modules $\mathcal O\subset \mathcal O+\mathcal Oa\subset\kappa$. Tensoring with $\reallywidehat{\mathcal O}$ over $\mathcal O$, we get inclusions $\reallywidehat{\mathcal O}=\reallywidehat{\mathcal O}+\reallywidehat{\mathcal O}a\subset\reallywidehat{\mathcal O}\otimes_{\mathcal O} \kappa$. As the homomorphism $\mathcal O\rightarrow \reallywidehat{\mathcal O}$ is faithfully flat, it follows that the inclusion $\mathcal O\subset \mathcal O+\mathcal Oa$ is as well an identity and therefore we have $a\in \mathcal O$; thus $\kappa\cap\reallywidehat{\mathcal O}=\mathcal O$.

As $\mathcal M$ is a torsion free $\mathcal O$-module, the same argument shows that we have 
$$\mathcal M=\mathcal O\otimes_{\mathcal O} \mathcal M=(\kappa\otimes_{\mathcal O} \mathcal M)\cap (\reallywidehat{\mathcal O}\otimes_{\mathcal O} \mathcal M)=V\cap (\reallywidehat{\mathcal O}\otimes_{\mathcal O} \mathcal M).$$ But as $\xi\in H^0(\reallywidehat{\mathbb P^N_l},\reallywidehat{\mathcal V})\subset V$, we also have $\xi\in V\cap (\reallywidehat{\mathcal O}\otimes_{\mathcal O} \mathcal M)$ and thus we have $\xi\in \mathcal M$. As $z$ is an arbitrary closed point of $\mathcal C$, we get that there exists an open subscheme $\mathcal W_{\xi}$ of $\mathcal U$ which contains $\mathcal C$ and such that we have $\xi\in H^0(\mathcal W_{\xi},\mathcal V)$. 

In this paragraph we consider the case when $\mathcal V$ is a reflexive $\mathcal O_{\mathcal U}$-module. As $\mathcal C$ has positive dimension, the complement of $\mathcal W_{\xi}$ in $\mathbb P^N_l$ (as it does not intersect $\mathcal C$) has codimension in $\mathbb P^N_l$ at least $2$. This implies that $\xi\in H^0(\mathcal W_{\xi},\mathcal V)=H^0(\mathcal U,\mathcal V)$, cf. Fact \ref{F1} (c) for the equality part. Thus we can take $\mathcal W_{\xi}=\mathcal U$. This implies that $\rho_{\mathcal U}$ is surjective and therefore an isomorphism. It is well known that $H^0(\mathcal U,\mathcal V)$ is a finite dimensional $l$-vector space (this can be checked by considering a coherent $\mathcal O_{\mathbb P^N_l}$-module that extends $\mathcal V$ and by recalling that $\mathbb P^N_l\setminus \mathcal U$ has codimension in $\mathbb P^N_l$ at least $2$). We conclude that $H^0(\reallywidehat{\mathbb P^N_l},\reallywidehat{\mathcal V})$ is a finite dimensional $l$-vector space when $\mathcal V$ is a reflexive $\mathcal O_{\mathcal U}$-module. 

We now consider the general case when $\mathcal V$ is torsion free but not necessarily a reflexive $\mathcal O_{\mathcal U}$-module. As $\mathcal V$ is torsion free, it is a $\mathcal O_{\mathcal U}$-submodule of a coherent locally free $\mathcal O_{\mathcal U}$-module. From this and the previous paragraph we get that $H^0(\reallywidehat{\mathbb P^N_l},\reallywidehat{\mathcal V})$ is a finite dimensional $l$-vector space. Let $\xi_1,\ldots, \xi_n\in H^0(\reallywidehat{\mathbb P^N_l},\reallywidehat{\mathcal V})\subset V$ be such that they generate $H^0(\reallywidehat{\mathbb P^N_l},\reallywidehat{\mathcal V})$ as an $l$-vector space. For each $i\in\{1,\ldots,n\}$ let $\mathcal W_{\xi_i}$ be an open subscheme of $\mathcal U$ which contains $\mathcal C$ and such that we have $\xi_i\in H^0(\mathcal W_{\xi_i},\mathcal V)$ (see above). Then for the open subscheme $\mathcal W:=\cap_{i=1}^n \mathcal W_i$ of $\mathcal U$, $\rho_{\mathcal W}$ is surjective and therefore also an isomorphism between finite dimensional $l$-vector spaces.
\endproof 

\section{Proof of Theorem \ref{T2}}\label{S5} 

In this section we prove Theorem \ref{T2}. By enlarging $l^{\prime}$ and making $\mathcal C$ smaller if needed, to prove Theorem \ref{T2} we can assume that $\mathcal C$ is a geometrically integral curve whose normalization $\mathcal C^{\textup{n}}$ is smooth over $\Spec(l^{\prime})$. Let $\mathcal I$ be the ideal sheaf of $\mathcal O_{\mathbb P^N_{l^{\prime}}}$ which defines $\mathcal C$. We have the following general constancy result. 

\begin{prop}\label{P6}
In the context of Theorem \ref{T2}, we assume that $l^{\prime}=l$ and that $\mathcal C$ is a geometrically integral curve with smooth normalization. Then the $\BT$ $\mathcal D_{\reallywidehat{\mathbb P^N_l}}$ over the formal completion $\reallywidehat{\mathbb P^N_l}$ of $\mathbb P^N_l$ along $\mathcal C$ induced naturally by the inductive system $\mathcal D_{n,\mathcal U_{\infty}}$ is as well constant, i.e., it is isomorphic to
$G_{\reallywidehat{\mathbb P^N_l}}$ with $G$ a $\BT$ over $\Spec(l)$.
\end{prop}

\noindent
{\bf Proof:} We first consider the case when $\mathcal C=\mathcal C^{\textup{n}}$ is smooth over $\Spec(l)$. The conormal sheaf $\mathcal I/\mathcal I^2$ of $\mathcal C$ inside $\mathbb P^N_l$ is a subsheaf of $\mathcal O_{\mathcal C}\otimes_{\mathcal O_{\mathbb P^N_l}}\Omega^1_{\mathbb P^N_l}$ and thus (cf. Euler's sequence) also of $\mathcal O_{\mathcal C}(-1)^{N+1}$. Therefore for each $t\in\mathbb N^*$, its $t$-th symmetric power $S^t(\mathcal I/\mathcal I^2)=\mathcal I^t/\mathcal I^{t+1}$ has no non-zero global section. 

We check by induction on $t\in\mathbb N$ that the $\BT$ given by the inductive system $\mathcal D_{n,\mathcal C_t}$ over the $t$-th infinitesimal neighborhood $\mathcal C_t$ of $\mathcal C$ in $\mathbb P^N_l$ is constant. We recall the convention of \cite{G4}, Def. 16.1.2 that $\mathcal C_t$ is defined by the ideal $\mathcal I^{t+1}$ of $\mathbb P^N_l$. As the inductive system $\mathcal D_{n,\mathcal C}$ is a constant $\BT$ over $\mathcal C=\mathcal C_0$ isomorphic to $G_{\mathcal C}$ for some $\BT$ $G$ over $\Spec(l)$, the base of the induction holds for $t=0$. As $G_{\mathcal C_{t+1}}$ is a canonical lift of $G_{\mathcal C_t}$ to $\mathcal C_{t+1}$, the lifts of $G_{\mathcal C_t}$ to $\BT$'s over $\mathcal C_{t+1}$ are parametrized by the global sections of an $\mathcal O_{\mathcal C}$-module isomorphic to $(\mathcal I^{t+1}/\mathcal I^{t+2})^{e_G}$, where $e_G$ is the product of the dimension and the codimension of $G$. From this and the previous paragraph we get that $G_{\mathcal C_{t+1}}$ is the only lift of $G_{\mathcal C_t}$ to a $\BT$ over $\mathcal C_{t+1}$. Thus if the inductive system $\mathcal D_{n,\mathcal C_t}$ is a constant $\BT$ over $\mathcal C_t$ isomorphic to $G_{\mathcal C_t}$, then the inductive system $\mathcal D_{n,\mathcal C_{t+1}}$ is as well a constant $\BT$ over $\mathcal C_{t+1}$ isomorphic to $G_{\mathcal C_{t+1}}$. This ends the inductive step and thus also the induction.

From the previous paragraph we get that the $\BT$ $\mathcal
D_{\reallywidehat{\mathbb P^N_l}}$ over the formal completion $\reallywidehat{\mathbb P^N_l}$ of $\mathbb
P^N_l$ along $\mathcal C$ induced naturally by the inductive system $\mathcal
D_{n,\mathcal U_{\infty}}$ is as well constant isomorphic to
$G_{\reallywidehat{\mathbb P^N_l}}$.

We now consider the case when $\mathcal C\neq\mathcal C^{\textup{n}}$. Let $\varepsilon:\mathcal C^{\textup{n}}\rightarrow \mathbb P^N_l\times_{\Spec(l)} \mathbb P^3_l$ be a closed embedding whose projections are the composite morphism $\mathcal C^{\textup{n}}\rightarrow\mathcal C\rightarrow \mathbb P^N_l$ and a closed embedding $\mathcal C^{\textup{n}}\rightarrow \mathbb P^3_l$. Let $\mathcal J$ be the ideal sheaf of
$\mathcal O_{\mathbb P^N_l\times_{\Spec(l)} \mathbb P^3_l}$ which defines $\mathcal C^{\textup{n}}$. We have a commutative diagram
\[\xymatrixcolsep{8pc}\xymatrix{
\mathcal C^{\textup{n}} \ar[r]^{\varepsilon} \ar[d]^{\textup{nat}} & \mathbb P^N_l\times_{\Spec(l)} \mathbb P^3_l \ar[d]^{1^{\textup{st}}\;\textup{projection}} \\
\mathcal C \ar[r]^{\textup{incl}} & \mathbb P^N_l. }\]
It is easy to see that for each $t\in \mathbb N^*$, the $t$-th symbolic power $\mathcal I^{(t)}$ of $\mathcal I$ (in the sense of \cite{ZS1}, Ch. IV, Def. of Sect. 12) is the largest $\mathcal O_{\mathbb P^N_l}$-submodule of $\mathcal I$ that maps naturally into $\mathcal J^t$. As above we argue that the $t$-th symmetric power $S^t(\mathcal J/\mathcal J^2)=\mathcal J^t/\mathcal J^{t+1}$ has no non-zero global section. From this and the natural inclusion $\mathcal I^{(t)}/\mathcal I^{(t+1)}\rightarrow \textup{nat}_*(\mathcal J^t/\mathcal J^{t+1})$ of $\mathcal O_{\mathcal C}$-modules, we get that also $\mathcal I^{(t)}/\mathcal I^{(t+1)}$ has no non-zero global section.\footnote{The fact that $H^0(\mathcal C, \mathcal J^{(t)}/\mathcal J^{(t+1)})=0$ for $t>0$ can be proved without using $\mathcal C^{\textup{n}}$ as follows. If $\mathcal Z$ is an integral closed scheme of a smooth $k$-scheme $\mathcal X$, $\mathcal J$ the ideal sheaf defining $\mathcal Z$ in $\mathcal X$, $\mathcal J^{(t)}$ the $t$-th symbolic power, $\textup{Diff}$ the sheaf  of differential operators on $\mathcal X$ in the sense of \cite{G4}, Sect. 16.8, filtered by $\textup{Diff}^n=\textup{Diff}^n_{\mathcal X/\Spec(k)}$ ($n\in\mathbb N$), then $\textup{Diff}^n(\mathcal J^{(t)})\subset \mathcal J^{(t-n)}$ ($\mathcal J^{(m)}:=\mathcal O_{\mathcal X}$ for $m\le 0$) and this induces $\textup{Diff}^t/\textup{Diff}^{t-1}\otimes_{\mathcal O_{\mathcal X}} \mathcal J^{(t)}/\mathcal J^{{(t+1)}}=\underline{Hom}(\textup{Sym}^t(\Omega^1_{\mathcal X/k}),\mathcal O_{\mathcal X})\rightarrow\mathcal O_{\mathcal Z}$, hence $\mathcal J^{(t)}/\mathcal J^{{(t+1)}}\rightarrow\textup{Sym}^t(\Omega^1_{\mathcal X/k})\otimes_{\mathcal O_{\mathcal X}} \mathcal O_{\mathcal Z}$ which is injective if $\mathcal Z$ is generically smooth. For $\mathcal X$ a projective space and $\mathcal Z$ of positive dimension and generically smooth, this proves that $H^0(\mathcal Z, \mathcal J^{(t)}/\mathcal J^{(t+1)})=0$ for $t>0$.} 

As the inductive system $\mathcal D_{n,\mathcal C}$ is a constant $\BT$ over $\mathcal C$ isomorphic to $G_{\mathcal C}$, an induction on $t\in\mathbb N$ similar to the one above shows that the $\BT$ given by the inductive system $\mathcal D_{n,\mathcal C_{(t)}}$ over the $t$-th symbolic infinitesimal neighborhood $\mathcal C_{(t)}$ of $\mathcal C$ in $\mathbb P^N_l$ (defined by the ideal sheaf $\mathcal I^{(t+1)}$ of $\mathbb P^N_l$), is a constant $\BT$ over $\mathcal C_{(t)}$ isomorphic to $G_{\mathcal C_{(t)}}$, compatibly for different $t$'s. 

By \cite{ZS2}, Ch. VIII, Cor. 5 of Thm. 13 applied to inclusions $\mathcal I(U)\subset\mathcal O_{\mathbb P^N_l}(U)$ for finitely many affine open subschemes $U$ of $\mathbb P^N_l$, we have that for each $t\in\mathbb N^*$ there exists an integer $q\ge t$ such that we have $\mathcal I^{(q)}\subset\mathcal I^t$. This implies that 
$$\mathop{\lim_{\longleftarrow}}_{t\in\mathbb N^*}\mathcal O_{\mathbb P^N_l}/\mathcal I^t=\mathop{\lim_{\longleftarrow}}_{t\in\mathbb N^*}\mathcal O_{\mathbb P^N_l}/\mathcal I^{(t)}.$$ 
From this and the previous paragraph we get that the $\BT$ $\mathcal D_{\reallywidehat{\mathbb P^N_l}}$ over the formal completion $\reallywidehat{\mathbb P^N_l}$ of $\mathbb P^N_l$ along $\mathcal C$ induced naturally by the inductive system $\mathcal D_{n,\mathcal U_{\infty}}$ is as well constant and isomorphic to $G_{\reallywidehat{\mathbb P^N_l}}$.\endproof

\medskip
To continue the proof of Theorem \ref{T2}, in this paragraph we consider the particular case when $l^{\prime}=l$ is an algebraically closed field. This assumption implies that $\mathcal C$ is G3 in $\mathbb P^N_l$, cf. Section \ref{S4}. From this and Lemma \ref{L7} applied with $(\mathcal U,\mathcal V)=(\mathcal U_n,\mathcal V_n)$, where $\mathcal V_n$ is the coherent locally free $\mathcal O_{\mathcal U_n}$-algebra that defines $\mathcal D_n$, we get that we have a natural ring identification $H^0(\mathcal U_n,\mathcal V_n)=H^0(\reallywidehat{\mathbb P^N_l},\reallywidehat{\mathcal V_n})$, where $\reallywidehat{\mathcal V_n}:=\chi^*(\mathcal V_n)$ is the formal completion of $\mathcal V_n$ along $\mathcal C$. From Proposition \ref{P6} we get that $H^0(\reallywidehat{\mathbb P^N_l},\reallywidehat{\mathcal V_n})$ is canonically identified with the ring of regular functions on $G[p^n]$. We have similar identifications for $\mathcal V_n\otimes_{\mathcal O_{\mathcal U_n}}\mathcal V_n$, $\reallywidehat{\mathcal V_n}\otimes_{\mathcal O_{\reallywidehat{\mathbb P^N_l}}}\reallywidehat{\mathcal V_n}$ and $G[p^n]\times_{\Spec(l)} G[p^n]$, with certain compatibilities. Therefore we have a canonical homomorphism $\vartheta_n:\mathcal D_n\rightarrow G[p^n]_{\mathcal U_n}$ of finite flat group schemes whose pullback to the formal completion $\reallywidehat{\mathbb P^N_l}$ is an isomorphism, being the truncation of level $n$ of an isomorphism as in Proposition \ref{P6}. From this and the first part of Fact \ref{F1} (d) we get that $\vartheta_n$ is an isomorphism over an open subscheme $\mathcal W_n$ of $\mathcal U_n$ which contains $\mathcal C$. As the complement of $\mathcal W_n$ in $\mathbb P^N_l$ has codimension in $\mathbb P^N_l$ at least $2$, from the second part of Fact \ref{F1} (d) we get that in fact $\vartheta_n$ is an isomorphism. The isomorphisms $\vartheta_n$'s are compatible in the sense that for all $n,m\in\mathbb N^*$ we have $\vartheta_{n+m}[p^n]=\vartheta_{n,\mathcal U_{n+m}}$. Therefore $\mathcal D_n$ extends uniquely (up to unique isomorphism) to a constant $\BT_n$ $\mathcal D^+_n$ over $\mathbb P^N_l$ isomorphic to $G[p^n]_{\mathbb P^N_l}$. Due to the mentioned compatibility of $\vartheta_n$'s and the uniqueness of the extension $\mathcal D^+_n$ of $\mathcal D_n$, we have a canonical identification $\mathcal D_{n+1}^+[p^n]=\mathcal D_n^+$. Thus the inductive system $\mathcal D_n^+$ is a constant $\BT$ over $\mathbb P^N_l$ isomorphic to $G_{\mathbb P^N_l}$. Moreover, $G[p^n]$ is the affine group scheme over $\Spec(l)$ defined by the Hopf $l$-algebra of global functions on $\mathcal D_n^+$ and thus also on either $\mathcal D_n$ or $\mathcal D_{n,\mathcal U_{n+m}}$, cf. Fact \ref{F1} (c) and the fact that the complement of $\mathcal U_n$ in $\mathbb P^N_l$ has codimension in $\mathbb P^N_l$ at least $2$. 

We now consider the general case. By enlarging $l^{\prime}$, we can assume that it is an algebraically closed field. Let $n,m\in\mathbb N^*$. If $\omega_n:\mathcal D_n\rightarrow \mathcal O_{\mathcal U_n}$ is the structure morphism, then the $\mathcal O_{U_n}$-linear map $\mathcal O_{\mathcal U_n}\otimes_l H^0(\mathcal U_n,\omega_{n,*}(\mathcal O_{\mathcal D_n}))\rightarrow \omega_{n,*}(\mathcal O_{\mathcal D_n})$ is an isomorphism (as the base change to $\mathcal U_{n,l^\prime}$ is) and thus the Hopf $\mathcal O_{U_n}$-algebra structure on $\omega_{n,*}(\mathcal O_{\mathcal D_n})$ defines a commutative and cocommutative Hopf $l$-algebra structure on $H^0(\mathcal U_n,\omega_{n,*}(\mathcal O_{\mathcal D_n}))$ and hence a finite group scheme $G_n$ over $l$ and a canonical isomorphism $\mathcal D_n\rightarrow G_{n,\mathcal U_n}$. This implies that $G_n$ is a $\BT_n$ over $\Spec(l)$. The canonical closed embedding homomorphisms $\mathcal D_{n,\mathcal U_{n+m}}\rightarrow \mathcal D_{n+m}$ and the canonical epimorphisms $\mathcal D_{n+m}\rightarrow \mathcal D_{n,\mathcal U_{n+m}}$ induce naturally homomorphisms $G_n\rightarrow G_{n+m}$ and $G_{n+m}\rightarrow G_n$ which are closed embeddings and epimorphisms (respectively) as their extensions to $l^{\prime}$ are so, cf. previous paragraph. These last homomorphisms define a $\BT$ $G$ over $\Spec(l)$ such that for all $n\in\mathbb N^*$ we have $G_n=G[p^n]$. We conclude that for each $n\in\mathbb N^*$, $\mathcal D_n$ extends uniquely (up to unique isomorphism) to a constant $\BT_n$ $\mathcal D_n^+$ over $\mathbb P^N_l$ isomorphic to $G[p^n]_{\mathbb P^N_l}$ and the identification $\mathcal D_{n+1}[p^n]=\mathcal D_{n,\mathcal U_{n+1}}$ extends to an identification $\mathcal D_{n+1}^+[p^n]=\mathcal D_n^+$ which is the pullback to $\mathbb P^N_l$ of the identification $(G[p^{n+1}])[p^n]=G[p^n]$ over $\Spec(l)$. Thus Theorem \ref{T2} holds.\endproof

\section{Proof of Theorem \ref{T1} (b)}\label{S6}

If $d=2$, then $R$ admits a faithfully flat regular local extension of the form $W(l)[[x_1,x_2]]/(h)$ with $l$ a perfect field which contains $k$ and such that $\bar h\notin (x_1^p,x_2^p,x_1^{p-1}x_2^{p-1})$, cf. Fact \ref{F4}. Thus $R$ is $p$-quasi-healthy by \cite{VZ2}, Thm. 3. 

In the rest of this section we will prove by induction on $d\ge 3$ that Theorem \ref{T1} (b) also holds for $d\ge 3$. The base of the induction (i.e., the case when $d=3$) is checked in Subsection \ref{SS61}. The inductive step (i.e., the passage from $d-1\ge 3$ to $d\ge 4$) is checked in Subsection \ref{SS62}. 

Let $D_U$ be a $\BT$ over the punctured spectrum $U$ of a regular local scheme $X=\Spec(R)$ of mixed characteristic $(0,p)$ and dimension $d\ge 3$. We have to show that if condition ($\natural$) holds for $R$, then $D_U$ extends to a $\BT$ $D$ over $X$. Based on Lemma \ref{L3} and Fact \ref{F4}, by passing to a faithfully flat extension we can assume that $R$ is also complete of the form $R=W(k)[[x_1,\ldots, x_d]]/(h)$, with $k$ an algebraically closed field of positive characteristic $p$ and with $h\in W(k)[[x_1,\ldots,x_d]]$ such that $\bar h\notin (x_1^p,\ldots,x_d^p)+(x_1,\ldots,x_d)^{2p-2}$. For $i\in\{1,\ldots,d\}$ let $y_i:=x_i+(h)\in R$. We can also assume that the residue field $k$ of $R$ is uncountable and (cf. Section \ref{S3}) that $y_1,\ldots,y_d$ is a regular system of parameters of $R$. 

\subsection{The base of the induction, i.e., the case $d=3$}\label{SS61} 
In this subsection we will assume that $d=3$. Thus $R=W(k)[[x_1,x_2,x_3]]/(h)$ has $y_1,y_2,y_3$ as a regular system of parameters of $R$. Let $Z$ be the blow up of $X$ along its closed point; it is a regular scheme of dimension $3$ which is projective over $X$ and which is the union of the open subscheme $U$ and of a closed subscheme $\mathbb P^2_k$. From Proposition \ref{P2} we get that the restriction of $D_U$ to the generic point of $Z$ (i.e., to $\Spec(\Frac(R))$) extends to a $\BT$ $D_{\Spec(O)}$ over the spectrum of the discrete valuation ring $O$ which is the local ring in $Z$ of the generic point of $\mathbb P^2_k$. From this and Fact \ref{F1} (b) we get that for each $n\in\mathbb N^*$ there exists a finite set $\mathcal S_n$ of closed points of $\mathbb P^2_k$ such that the $\BT_n$ $D_U[p^n]$ over $U$ extends uniquely (up to unique isomorphism) to a finite flat group scheme $E_n$ over $Z\setminus \mathcal S_n$ whose pullback to $\Spec(O)$ is $D_{\Spec(O)}[p^n]$. 

Recall the standard fact that if $\triangle$ is a finite locally free commutative group scheme annihilated by $p^n$ over a scheme $\Sigma$, then the set of points $z$ of $\Sigma$ such that the fiber of $\triangle$ at $z$ is a $\BT_n$ over the residue field of $z$ is open and $\triangle$ is a $\BT_n$ over the corresponding open subscheme. Thus, let $\mathcal K_n$ be the smallest reduced closed subscheme of $\mathbb P^2_k$ which contains $\mathcal S_n$ and such that the restriction of $E_n$ to $Z\setminus\mathcal K_n$ is a $\BT_n$; the dimension of $\mathcal K_n$ is at most $1$. 

We can assume we have a chain of inclusions $\mathcal S_1\subset\mathcal S_2\subset \cdots\subset\mathcal S_n\subset\cdots $. Thus we have a second chain of inclusions $\mathcal K_1\subset\mathcal K_2\subset \cdots\subset\mathcal K_n\subset\cdots $ such that for each $n\in\mathbb N^*$ the $\BT_n$'s $E_{n+1}[p^n]_{Z\setminus\mathcal K_{n+1}}$ and $E_{n,Z\setminus\mathcal K_{n+1}}$ coincide over $Z\setminus\mathcal K_{n+1}$. The set $\mathcal S_{\infty}:=\cup_{n\ge 1} \mathcal S_n$ is countable and the ind-constructible set $\mathcal K_{\infty}:=\cup_{n\ge 1} \mathcal K_n$ has a countable number of maximal points, and is the union of their closures.\footnote{Recall that a point of a topological space is said to be a maximal point if its closure is not strictly contained in the closure of another point. For sober spaces the maximal points are the generic points of the irreducible components.} 
Moreover, $D_U$ extends to a $\BT$ $D_{Z\setminus\mathcal K_{\infty}}$ over the stable under generization, pro-constructible subset $Z\setminus\mathcal K_{\infty}=\cap_{n=1}^{\infty} (Z\setminus\mathcal K_n)$ of $Z$ (see footnote 2). 

\begin{claim}\label{CL1}
For each $n\in\mathbb N^*$, $\mathcal K_n$ is a finite set and therefore (by enlarging $\mathcal S_n$) we can assume that $\mathcal S_n=\mathcal K_n$.
\end{claim}

We begin the proof of Claim \ref{CL1} by introducing notation pertaining to $h$.

\subsubsection{On $\bar h\in k[[x_1,x_2,x_3]]$}\label{SSS611}
For each integer $i\in\{1,\ldots,2p-3\}$ let $\mathbb J_i^3$ be the subset of $\{0,\ldots,p-1\}^3$ formed by those triples whose sum is $i$. We write 
$$\bar h(x_1,x_2,x_3)=\bar h_0(x_1,x_2,x_3)+\sum_{i=1}^{2p-3} \bar h_i(x_1,x_2,x_3),$$ 
where $\bar h_0$ belongs to the ideal $(x_1^p,x_2^p,x_3^p)+(x_1,x_2,x_3)^{2p-2}$ of $k[[x_1,x_2,x_3]]$ and where each 
$$\bar h_i(x_1,x_2,x_3)=\sum_{(i_1,i_2,i_3)\in\mathbb J_i^3} \delta_{i_1,i_2,i_3}x_1^{i_1}x_2^{i_2}x_3^{i_3}\in k[x_1,x_2,x_3]$$ 
 is either $0$ or a homogeneous polynomial of degree $i$. As $\bar h$ does not belong to the ideal $(x_1^p,x_2^p,x_3^p)+(x_1,x_2,x_3)^{2p-2}$ of $k[[x_1,x_2,x_3]]$, there exists a smallest element $c\in\{1,\ldots,2p-3\}$ such that we have $\bar h_c\neq 0$. 

\subsubsection{Good regular closed subschemes of $X$ of dimension $2$}\label{SSS62}
For each triple $\zeta=(\zeta_1,\zeta_2,\zeta_3)\in W(k)^3$ such that its reduction $\bar\zeta=(\bar\zeta_1,\bar\zeta_2,\bar\zeta_3)$ modulo $p$ is a non-zero element of $k^3$, we consider the regular closed subscheme $\Spec(S_{\zeta})$ of $X$, where $S_{\zeta}:=R/(\zeta_1y_1+\zeta_2y_2+\zeta_3y_3)$. Let $Z_{\zeta}$ be the
closed subscheme of $Z$ which is the blow up of $\Spec(S_{\zeta})$ along its closed point; it is the union of its open
subscheme $U\cap\Spec(S_{\zeta})$ and of a $\mathbb P^1_{k,\zeta}$ curve inside $\mathbb
P^2_k$. It is easy to see that, with respect to the usual projective coordinates $[w_0,w_1,w_2]$ of $\mathbb P^2_k$, the curve $\mathbb P^1_{k,\zeta}$ of $\mathbb P^2_k$ is defined by the equation $\bar\zeta_1w_0+\bar\zeta_2w_1+\bar\zeta_3w_2=0$. If $P=[\gamma_0,\gamma_1,\gamma_2]\in \mathcal S_{\infty}\subset \mathbb P^2_k$ with $\gamma_0,\gamma_1,\gamma_2\in k$ not all zero, then the condition $P\notin \mathbb P^1_{k,\zeta}$ gets translated into the inequality $\bar\zeta_1\gamma_0+\bar\zeta_2\gamma_1+\bar\zeta_3\gamma_2\neq 0$. Thus there exists a countable union $\mathcal L_{\infty}$ of lines in $\mathbb P^2_k$ such that we have $\mathbb P^1_{k,\zeta}\cap \mathcal S_{\infty}=\emptyset$ if and only if $[\bar\zeta_1,\bar\zeta_2,\bar\zeta_3]\notin \mathcal L_{\infty}(k)$. 

Let $\mathcal P_{\infty}$ be the countable subset of $\mathbb P^2_k(k)$ such that $\mathbb P^1_{k,\zeta}$ is not contained in $\mathcal K_{\infty}$ if and only if $[\bar\zeta_1,\bar\zeta_2,\bar\zeta_3]\notin \mathcal P_{\infty}$.

From Lemma \ref{L6} applied with $(l,d)=(k,3)$, we get that there exists an open dense subscheme $\mathbb O$ of $\mathbb A^2_k$ such that for each pair $(v_1,v_2)\in \mathbb O(k)\subset k^2$ the polynomial
$$\bar h_{c,\zeta}(x_1,x_2)=\bar h_c(x_1,x_2,v_1x_1+v_2x_2)$$ 
does not belong to the ideal $(x_1^p,x_2^p,x_1^{p-1}x_2^{p-1})$ of $k[[x_1,x_2]]$. If $(\bar\zeta_1,\bar\zeta_2,\bar\zeta_3)=(-v_1,-v_2,1)$ with $(v_1,v_2)\in \mathbb O(k)$, then the regular ring $S_{\zeta}=R/(\zeta_1y_1+\zeta_2y_2+\zeta_3y_3)$ is isomorphic to $W(k)[[x_1,x_2]]/(h_{\zeta})$, where 
$$h_{\zeta}(x_1,x_2)=h(x_1,x_2,-\zeta_3^{-1}\zeta_1x_1-\zeta_3^{-1}\zeta_2x_2)\in W(k)[[x_1,x_2]]$$ 
is such that its reduction $\bar h_{\zeta}$ modulo $p$ has $\bar h_{c,\zeta}$ as its homogeneous component of degree $c$ and therefore it does not belong to the ideal $(x_1^p,x_2^p,x_1^{p-1}x_2^{p-1})$ of $k[[x_1,x_2]]$; thus in such a case $S_{\zeta}$ is of mixed characteristic $(0,p)$ and is $p$-quasi-healthy (by \cite{VZ2}, Thm. 3).

We identify $\mathbb A^2_k$ with an open subscheme of $\mathbb P^2_k$ via the embedding $(v_1,v_2)\rightarrow [-v_1,-v_2,1]$. As $k$ is uncountable, the set of all closed points of the open dense subscheme $\mathbb O$ of $\mathbb A^2_k\subset\mathbb P^2_k$ cannot be contained in $\mathcal L_{\infty}(k)\cup\mathcal P_{\infty}$, i.e., there exist pairs $(v_1,v_2)\in\mathbb O(k)$ such that $[-v_1,-v_2,1]\notin \mathcal L_{\infty}(k)\cup\mathcal P_{\infty}$. If $(v_1,v_2)$ is such a pair, then by choosing $(\bar\zeta_1,\bar\zeta_2,\bar\zeta_3)=(-v_1,-v_2,1)$ we get that the regular closed subscheme $\Spec(S_{\zeta})$ of $X$ is good in the following sense:

\medskip
($\sharp$) {\it The regular ring $S_{\zeta}$ is $p$-quasi-healthy and moreover the line $\mathbb P^1_{k,\zeta}$ neither intersects the countable subset $\mathcal S_{\infty}$ of $\mathbb P^2_k$ nor is contained in $\mathcal K_{\infty}$.} 

\subsubsection{Proof of the Claim \ref{CL1}}\label{SSS63}
As $S_{\zeta}$ is $p$-quasi-healthy (cf. ($\sharp$)), the restriction of $D_U$ to $U\cap\Spec(S_{\zeta})$
extends uniquely (up to unique isomorphism) to a $\BT$ $D_{\Spec(S_{\zeta}})$ over $\Spec(S_{\zeta})$. As $\mathbb P^1_{k,\zeta}$ does not intersect $\mathcal S_n$ (cf. ($\sharp$)), for $n\in \mathbb N^*$ we can speak about the pullback $E_{n,\zeta}$ of $E_n$ to $Z_{\zeta}$. Let $E_{n,O_{\zeta}}$ be the restriction of $E_{n,\zeta}$ to the spectrum of the local ring $O_{\zeta}$ of $Z_{\zeta}$ which is a discrete valuation ring that dominates $S_{\zeta}$. As $\mathbb P^1_{k,\zeta}$ is not contained in $\mathcal K_{\infty}$, the inductive system $E_{n,O_{\zeta}}$ is a $\BT$ over $\Spec(O_{\zeta})$ which, based on Tate's extension theorem, is the pullback of $D_{\Spec(S_{\zeta})}$ to $\Spec(O_{\zeta})$. This implies that for each $n\in \mathbb N^*$, $E_{n,\zeta}$ is the pullback to $Z_{\zeta}$ of $D_{\Spec(S_{\zeta})}[p^n]$. Thus the inductive system $E_{n,\zeta}$ is a $\BT$ over $Z_{\zeta}$. This implies that $Z_{\zeta}$ and therefore also $\mathbb P^1_{k,\zeta}$ does not intersect $\mathcal K_{\infty}$. As two irreducible projective curves in $\mathbb P^2_k$ always intersect and as $\mathbb P^1_{k,\zeta}\cap\mathcal K_n=\emptyset$, we conclude that each $\mathcal K_n$ has dimension $0$ and therefore (by enlarging $\mathcal S_n$) we can assume that for each $n\in \mathbb N^*$ we have $\mathcal S_n=\mathcal K_n$. Thus also $\mathcal K_{\infty}=\mathcal S_{\infty}$ and each $E_n$ is a $\BT_n$ over $Z\setminus\mathcal S_n$. This ends the proof of the Claim \ref{CL1}.\endproof

\subsubsection{Applying Theorem \ref{T2}}\label{SSS64} 
Let $G$ be the fiber of $D_{\Spec(S_{\zeta})}$ over the closed point $\Spec(k)$ of $\Spec(S_{\zeta})$. For $n\in\mathbb N^*$, as $E_{n,\zeta}$ is the pullback of $D_{\Spec(S_{\zeta})}[p^n]$ to $Z_{\zeta}$, the pullback of $E_n$ to $\mathbb P^1_{k,\zeta}$ is canonically identified with $G[p^n]_{\mathbb P^1_{k,\zeta}}$. From this and Theorem \ref{T2} applied with
$(l,N,\mathcal C,\mathcal U_n,\mathcal D_n)=(k,2,\mathbb P^1_{k,\zeta},\mathbb
P^2_k\setminus\mathcal K_n,E_{n,\mathbb P^2_k\setminus\mathcal K_n})$ we get that the restriction $D_{\mathbb P^2_k\setminus\mathcal S_{\infty}}$ of $D_{Z\setminus\mathcal S_{\infty}}$ to $\mathbb P^2_k\setminus\mathcal S_{\infty}$ is isomorphic to $G_{\mathbb P^2_k\setminus\mathcal S_{\infty}}$ and thus it extends to a constant $\BT$ $D_{\mathbb P^2_k}$ over $\mathbb P^2_k$ isomorphic to $G_{\mathbb P^2_k}$.

\subsubsection{Liftings to infinitesimal neighborhoods of $\mathbb P^2_k$ in $Z$}\label{SSS65}
Let $\mathfrak m:=(y_1,y_2,y_3)$ be the maximal ideal of $R$. As $\mathfrak m\mathcal O_Z$ is the ideal sheaf of $\mathcal O_Z$ that defines the
closed subscheme $\mathbb P^2_k$ of $Z$, for each $t\in \mathbb N$ for the invertible $\mathcal O_{\mathbb P^2_k}$-module
$\mathfrak m^t\mathcal O_Z/\mathfrak m^{t+1}\mathcal O_Z$ we have 
$$H^1(Z,\mathfrak m^t\mathcal O_Z/\mathfrak m^{t+1}\mathcal
O_Z)=H^1(\mathbb P^2_k,\mathfrak m^t\mathcal O_Z/\mathfrak m^{t+1}\mathcal O_Z)=0.$$ 

We consider a coherent locally free module $\mathcal V_{t+1}$ of rank $r$ over the $t+1$-th infinitesimal neighborhood $\mathbb P^2_{k,t+1}$ of $\mathbb P^2_k$ in $Z$ (i.e., over the reduction modulo $\mathfrak m^{t+2}\mathcal O_Z$ of $Z$) such that the coherent locally free module
$\mathcal V_t:=\mathcal V_{t+1}/\mathfrak m^{t+1}\mathcal V_{t+1}$ over the $t$-th
infinitesimal neighborhood $\mathbb P^2_{k,t}$ of $\mathbb P^2_k$ in $Z$ is trivial and
thus isomorphic to $\mathcal O_{\mathbb P^2_{k,t}}^r$. As we have $\mathfrak m^{t+1}\mathcal V_{t+1}=(\mathfrak m^{t+1}\mathcal O_Z/\mathfrak m^{t+2}\mathcal O_Z)\otimes_{\mathcal O_Z} \mathcal V_{t+1}=(\mathfrak m^{t+1}\mathcal O_Z/\mathfrak m^{t+2}\mathcal O_Z)\otimes_{\mathcal O_Z} \mathcal V_t=(\mathfrak m^{t+1}\mathcal O_Z/\mathfrak m^{t+2}\mathcal O_Z)^r$, the short exact sequence 
$$0\rightarrow\mathfrak m^{t+1}\mathcal V_{t+1}\rightarrow \mathcal V_{t+1}\rightarrow\mathcal V_t\rightarrow 0$$ 
of coherent $\mathcal O_Z$-modules induces an exact complex 
$$H^0(Z,\mathcal V_{t+1})\rightarrow H^0(Z,\mathcal V_t)\rightarrow H^1(Z,\mathfrak m^{t+1}\mathcal O_Z/\mathfrak m^{t+2}\mathcal
O_Z)^r=0.$$ 
Therefore the reduction homomorphism $H^0(Z,\mathcal V_{t+1})\rightarrow H^0(Z,\mathcal V_t)$ is surjective and this implies that $\mathcal V_{t+1}$ is as well a free $\mathcal O_{\mathbb P^2_{k,t+1}}$-module of rank $r$.

By induction on $t\in \mathbb N$ we check that the constant $\BT$ $D_{\mathbb P^2_k}$ over $\mathbb P^2_k$ lifts uniquely (up to unique isomorphism) to a $\BT$ $D_{\mathbb P^2_{k,t}}$ over $\mathbb P^2_{k,t}$ in such a way that its restriction to the reduction modulo
$\mathfrak m^{t+1}$ of $Z\setminus \mathcal S_{\infty}$ is induced naturally by
$D_{Z\setminus\mathcal S_{\infty}}$ and moreover the coherent $\mathcal O_{\mathbb P^2_{k,t}}$-module associated naturally to the structure sheaf of the $\BT_n$ $D_{\mathbb P^2_{k,t}}[p^n]$ is free for all $n\in\mathbb N^*$ (in particular, each $D_{\mathbb P^2_{k,t}}[p^n]$ extends the reduction of $E_n$ modulo $\mathfrak m^{t+1}$).

The case $t=0$ was accomplished in Subsubsection \ref{SSS64}. The passage from $t$ to $t+1$ goes as follows. From the Grothendieck--Messing deformation theory we get that the lifts of the $\BT$ $D_{\mathbb P^2_{k,t}}$ to a $\BT$ over $\mathbb P^2_{k,t+1}$ are parametrized by the global sections of a torsor under the group of global sections of a coherent locally free $\mathcal O_{\mathbb P^2_k}$-module $\mathcal F_t$. Similarly, for each integer $\tilde t\ge 1$ the lifts of the $\BT_{t+\tilde t}$ which is the reduction modulo $\mathfrak m^{t+1}$ of $E_{t+\tilde t}$ to $\BT_{t+\tilde t}$'s over the reduction modulo $\mathfrak m^{t+2}$ of $Z\setminus \mathcal S_{t+\tilde t}$, are parametrized by a torsor under the group $H^0(\mathbb P^2_k\setminus \mathcal S_{t+\tilde t},\mathcal F_t)$ in a way compatible in $\tilde t$ (cf. \cite{I}, Thm. 4.4 c) and Cor. 4.7). As $H^0(\mathbb P^2_k\setminus \mathcal S_{t+\tilde t},\mathcal F_t)=H^0(\mathbb P^2_k\setminus \mathcal S_{\infty},\mathcal F_t)$ does not depend on $\tilde t$, we get that the lifts of the $\BT$ which is the restriction of $D_{Z\setminus\mathcal S_{\infty}}$ to the reduction modulo $\mathfrak m^{t+1}$ of $Z\setminus \mathcal S_{\infty}$ to $\BT$'s over the reduction modulo $\mathfrak m^{t+2}$ of $Z\setminus \mathcal S_{\infty}$, are parametrized by a torsor under the group $H^0(\mathbb P^2_k\setminus \mathcal S_{\infty},\mathcal F_t)$. As we have $H^0(\mathbb P^2_k,\mathcal F_t)=H^0(\mathbb P^2_k\setminus \mathcal S_{\infty},\mathcal F_t)$, we conclude that there exists a unique (up to unique isomorphism) $\BT$ over $\mathbb P^2_{k,t+1}$ which lifts $D_{\mathbb P^2_{k,t}}$ in such a way that its restriction to the reduction modulo $\mathfrak m^{t+2}$ of $Z\setminus \mathcal S_{\infty}$ is induced naturally by $D_{Z\setminus\mathcal S_{\infty}}$. Based on the previous paragraph, we get that the coherent $\mathcal O_{\mathbb P^2_{k,t+1}}$-module associated naturally to the structure sheaf of the $\BT_n$ $D_{\mathbb P^2_{k,t+1}}[p^n]$ is free for all $n\in\mathbb N^*$. This ends the induction on $t\in\mathbb N$. 

\subsubsection{End of the proof in the case $d=3$}\label{SSS66}
For $z\in\mathcal S_n$, let $\tau\in\mathcal O_{Z,z}$ be such that we have $\Spec(\mathcal O_{Z,z}/(\tau))=\Spec(\mathcal O_{Z,z})\times_Z \mathbb P^2_k$. As $\mathcal O_{Z,z}$ is a regular local ring of dimension $d=3$, for $t\in\mathbb N$ we have $\depth_{\mathcal O_{Z,z}} (\mathcal O_{Z,z}/\tau^{t+1}\mathcal O_{Z,z})=d-1\ge 2$ (cf. Auslander--Buchsbaum formula) and thus for each coherent locally free $\mathcal O_{\mathbb P^2_{k,t}}$-module $\mathfrak F_t$ we have a canonical identification $H^0(\mathbb P^2_{k,t}\setminus \mathcal S_n,\mathfrak F_t)=H^0(\mathbb P^2_{k,t},\mathfrak F_t)$. Thus the direct image via the open embedding $\mathbb P^2_{k,t}\setminus \mathcal S_n\rightarrow \mathbb P^2_{k,t}$ of the $\mathcal O_{\mathbb P^2_{k,t}\setminus \mathcal S_n}$-module associated to the reduction of $E_n$ modulo $\mathfrak m^{t+1}$ is the $\mathcal O_{\mathbb P^2_{k,t}}$-module associated to $D_{\mathbb P^2_{k,t}}[p^n]$. Based on this, from Lemma \ref{L4} applied to the local rings $\mathcal O_{Z,z}$ we get that for each $n\in\mathbb N^*$ the locally free $\mathcal O_{Z\setminus\mathcal S_n}$-module associated to $E_n$ extends (uniquely up to unique isomorphism) to a locally free $\mathcal O_Z$-module whose reduction modulo each $\mathfrak m^{t+1}$ is the $\mathcal O_{\mathbb P^2_{k,t}}$-module associated to $D_{\mathbb P^2_{k,t}}[p^n]$. This implies that each $E_n$ extends to a $\BT_n$ $E_n^+$ over $Z$ which lifts each $D_{\mathbb P^2_{k,t}}[p^n]$. For all $n,m\in\mathbb N^*$, the closed embedding homomorphisms $E_{n,Z\setminus\mathcal S_{n+m}}\rightarrow E_{n+m}$ extend to homomorphisms $E_n^+\rightarrow E_{n+m}^+$ over $Z$ which are closed embeddings identifying $E_n^+$ with $E_{n+m}^+[p^n]$, as their restrictions to $\mathbb P^2_{k,0}=\mathbb P^2_k$ are so. Thus the inductive system $E_n^+$ is a $\BT$ $E^+$ over $Z$. 

The $D_{\mathbb P^2_{k,t}}$'s define a $\BT$ $D_{\widehat{Z}}$ over the formal scheme $\widehat{Z}$ of the completion of $Z$ along $\mathbb P^2_k$ and moreover the coherent $\mathcal O_{\widehat{Z}}$-module $\mathcal O_{D_{\widehat{Z}}[p^n]}$ is free for all $n\in\mathbb N^*$. From \cite{G1}, Thm. 5.1.4 we get that $D_{\widehat{Z}}$ is the formal completion along $\mathbb P^2_k$ of a uniquely determined (up to unique isomorphism) $\BT$ $D_Z$ over $Z$ and that for all $n\in\mathbb N^*$ the coherent $\mathcal O_Z$-module $\mathcal O_{D_Z[p^n]}$ is free. Due to the uniqueness part we have $E^+=D_Z$. 

As the morphism $Z\rightarrow X$ is birational and projective and as $R$ is normal, the ring of global functions on $Z$ is $R$. From this and the fact that for all $n\in\mathbb N^*$ the coherent $\mathcal O_Z$-module $\mathcal O_{E_n^+}=\mathcal O_{D_Z[p^n]}$ is free, we get that $E^+=D_Z$ is the pullback of a uniquely determined (up to unique isomorphism) $\BT$ $D$ over $X$ which extends $D_U$ and whose truncations $D[p^n]$ are defined by identities $H^0(X,D[p^n])=H^0(Z,E^+[p^n])$ (to be compared with the Hopf algebra argument involving $\omega_n$'s at the end of Section \ref{S5}).

\subsection{The inductive step}\label{SS62}
In this subsection we will assume that $d\ge 4$ and that Theorem \ref{T1} (b) holds for regular rings of dimension at most $d-1$. As in Subsection \ref{SS61}, for each $d$-tuple $\zeta=(\zeta_1,\ldots,\zeta_d)\in W(k)^d$ such that its reduction modulo $p$ is not zero, we consider the regular ring $S_{\zeta}:=R/(\sum_{i=1}^d\zeta_iy_i)$. As $k$ is infinite, from Proposition \ref{P3} we get that we can choose $\zeta$ such that $S_{\zeta}$ is isomorphic to $W(k)[[x_1,x_2,\ldots,x_{d-1}]]/(h_{\zeta})$, where $h_{\zeta}\in W(k)[[x_1,x_2,\ldots,x_{d-1}]]$ is such that its reduction modulo $p$ does not belong to the ideal $(x_1^p,x_2^p,\ldots,x_{d-1}^p)+(x_1,x_2,\ldots,x_{d-1})^{2p-2}$ of $k[[x_1,x_2,\ldots,x_{d-1}]]$. We know that $S_{\zeta}$ is $p$-quasi-healthy, by the inductive assumption that Theorem \ref{T1} (b) holds for regular rings of dimension at most $d-1$. Thus $D_U$ modulo $(\sum_{i=1}^d\zeta_iy_i)$ extends to a $\BT$ over $\Spec(S_{\zeta})$. From this and Lemma \ref{L5} applied with $(R,y,R/(y),D_U)=(R,\sum_{i=1}^d\zeta_iy_i,S_{\zeta},D_U)$ we get that $D_U$ extends to a $\BT$ over $X$. Thus $R$ is $p$-quasi-healthy. This ends the induction and the proof of Theorem \ref{T1} (b).\endproof
 
\section{Proof of Theorem \ref{T1} (a) and Corollary \ref{C1}}\label{S7}
In the situation of Theorem \ref{T1}, from Theorem \ref{T1} (b) and Proposition \ref{P5} we get that each local ring of $X$ of mixed characteristic $(0,p)$ and of dimension at least $2$ is $p$-quasi-healthy. Thus Theorem \ref{T1} (a) is a particular case of the following general lemma whose proof relies on the purity of the branch locus.

\begin{lemma}\label{L8}
Let $Y$ be an integral scheme flat over $\mathbb Z$ such that $Y[\frac{1}{p}]$ is regular and each local ring of $Y$ is normal noetherian. We assume that the following two conditions hold:

\medskip\noindent
{\bf (i)} every local ring of $Y$ of mixed characteristic $(0,p)$ and dimension at least $2$ is $p$-quasi-healthy;

\smallskip\noindent
{\bf (ii)} there exists an affine open cover $(W_{\lambda})_{\lambda\in\Lambda}$ of $Y$ such that for each $\lambda\in\Lambda$ there exists $N_{\lambda}\in\mathbb N^*$ with the property that for every maximal point $\eta_0\in Y_{\mathbb F_p}\cap W_{\lambda}$, the absolute ramification index of $\mathcal O_{Y,\eta_0}$ is at most $N_\lambda$.% which is not a closed point

\medskip
If $D_{\eta}$ is a $\BT$ over the generic point $\eta$ of $Y$ which extends to every one dimensional local ring of $Y$, then $D_{\eta}$ extends to $Y$ (uniquely up to unique isomorphism).
\end{lemma}

\noindent
{\bf Proof:} As the local rings of $Y$ are normal noetherian, the same arguments as in the proof of \cite{T}, Thm. 4 give that the functor
$$(\textup{BT}\;\textup{groups}\;\textup{over}\; Y)\longrightarrow (\textup{BT}\;\textup{groups}\;\textup{over}\; \eta)$$
is fully faithful. Thus it is enough to prove the assertion locally and we can assume that $Y=\Spec(A)=W_{\lambda}$ (so $\Lambda=\{\lambda\}$). 

In this paragraph we check that the $\BT$ $D_{\eta}$ extends to a $\BT$ $D_{Y[\frac{1}{p}]}$ over $Y[\frac{1}{p}]=\Spec(A[\frac{1}{p}])$. For $n\in\mathbb N^*$, we consider the normalization $D_{Y[\frac{1}{p}],n}$ of $Y[\frac{1}{p}]$ in $D_{\eta}[p^n]$. As the finite \'etale group scheme $D_{\eta}[p^n]$ extends to a finite \'etale group scheme over the spectrum of each local ring of $Y[\frac{1}{p}]$ of dimension $1$ (i.e., of each local ring of $Y[\frac{1}{p}]$ which is a discrete valuation ring), the morphism $\vartheta_n:D_{Y[\frac{1}{p}],n}\rightarrow Y[\frac{1}{p}]$ is finite and \'etale over each local ring of $Y[\frac{1}{p}]$ of dimension $1$. From this and the purity of the branch locus (see \cite{G2}, Exp. X, Thm. 3.4 (i)]) we get that the morphism $\vartheta_n$ is finite and \'etale over each local ring of $Y[\frac{1}{p}]$. It is easy to see that this implies that $\vartheta_n$ defines a finite \'etale scheme over $Y[\frac{1}{p}]$ which extends $D_{\eta}[p^n]$ and thus has a unique group scheme structure which extends the group scheme structure on $D_{\eta}[p^n]$, and this defines a $\BT_n$ $D_{Y[\frac{1}{p}],n}$ over $Y[\frac{1}{p}]$. For $n,m\in\mathbb N^*$, the inclusions $D_{\eta}[p^n]\rightarrow D_{\eta}[p^{n+m}]$ extend to closed embedding homomorphisms $D_{Y[\frac{1}{p}],n}\rightarrow D_{Y[\frac{1}{p}],n+m}$ and the inductive system $D_{Y[\frac{1}{p}],n}$ is the unique (up to unique isomorphism) $\BT$ $D_{Y[\frac{1}{p}]}$ over $Y[\frac{1}{p}]$ which extends $D_{\eta}$.

To check that $D_{Y[\frac{1}{p}]}$ extends to $Y$ we consider two cases as follows.

\medskip
{\bf Case 1: $A/pA$ is noetherian.} Let $O_1,\ldots,O_t$ be all local rings of $R$ which are discrete valuation rings of mixed characteristic $(0,p)$; they correspond to the maximal points $\eta_1,\ldots,\eta_t$ (respectively) of $Y_{\mathbb F_p}$. For $i\in\{1,\ldots,t\}$ let $D_{O_i}$ be the unique (up to unique isomorphism) $\BT$ over $\Spec(O_i)$ which extends $D_{\eta}$. For each $n\in\mathbb N^*$ we consider the largest open subscheme $U_n$ of $Y$ which contains $Y[\frac{1}{p}]\cup \{\eta_1,\ldots,\eta_t\}$ and over which there exists a finite flat (locally free) commutative group scheme $D_{U_n,n}$ which extends compatibly $D_{Y[\frac{1}{p}]}[p^n]$ and $D_{O_i}[p^n]$ for all $i\in\{1,\ldots,t\}$.\footnote{In the general non-noetherian case ``finite flat" does not imply locally free, but in our case this follows from \cite{GR}, Prop. 2.4.19 applied to classical extensions of the type $A\subset A[\frac{1}{p}]$.} Let $W_n$ be the largest open subscheme of $U_n$ with the property that the restriction of $D_{U_n,n}$ to $W_n$ is a $\BT_n$ over $W_n$. We have a chain of inclusions $W_1\supset W_2\supset W_3\supset\cdots\supset W_m\supset\cdots$. Thus for $Y_i:=Y\setminus W_i$, we have a chain of inclusions $Y_1\subset Y_2\subset Y_3\subset\cdots\subset Y_m\subset\cdots$ between reduced closed subschemes of $Y_{\mathbb F_p}$. The codimension of each $Y_m$ is at least $2$. Note that the formation of the open subschemes $U_n$ and $W_n$ commutes with passage to spectra of local rings of $Y$.

We will show that the assumption that there exists $q\in\mathbb N^*$ such that $Y_q$ is non-empty leads to a contradiction. We can choose $q\in\mathbb N^*$ such that for all $m\in\mathbb N^*$ we have $\codim(Y_q)=\codim(Y_{q+m})=c\ge 2$. 

Let $z$ be a generic point of an irreducible component of $Y_q$ of codimension $c$ and let $\tilde R:=\mathcal O_{Y,z}$. We have $\dim(\tilde R)=c\ge 2$ and therefore $\tilde R$ is $p$-quasi-healthy. Let $\tilde X:=\Spec(\tilde R)$ and let $\tilde U$ be the punctured spectrum of $\tilde R$. From the very definitions and the choice of $z$ we get that $\tilde U\cap Y_n=\emptyset$ for all $n\in\mathbb N^*$. Thus there exists a unique (up to unique isomorphism) $\BT$ $D_{\tilde U}$ over $\tilde U$ which extends the restriction of $D_{Y[\frac{1}{p}]}$ to $\tilde X[\frac{1}{p}]=\Spec(\tilde R[{\frac{1}{p}}])$. Therefore $D_{\tilde U}$ extends uniquely (up to unique isomorphism) to a $\BT$ over $\tilde X$ and thus for each $n\in \mathbb N^*$ we have $z\notin Y_n$. This contradicts the fact that $z\in Y_q$. Thus our assumption leads to a contradiction. Therefore for all $n\in\mathbb N^*$ we have $U_n=W_n=Y$. This implies that $D_{\eta}$ extends to a $\BT$ over $Y$. 

\medskip
{\bf Case 2: general case.} Let $n\in\mathbb N^*$. Let $s\in\mathbb N$ be such that it depends only on $N_{\lambda}$ and \cite{B}, Thm. E or \cite{VZ3}, Cor. 3 applies to all discrete valuation rings which are local rings of $Y$ at maximal points of $Y_{\mathbb F_p}$ (see \cite{VZ3}, Exs. 2 and 4). Let $z$ be a point of $Y_{\mathbb F_p}$. From Case 1 applied to $\mathcal O_{Y,z}$ we get that the restriction of $D_{Y[\frac{1}{p}]}$ to $\Spec(\mathcal O_{Y,z}[\frac{1}{p}])$ extends to a $\BT$ $D_z$ over $\Spec(\mathcal O_{Y,z})$. Thus $D_z[p^{n+s}]$ extends to a $\BT_{n+s}$ $E_{n+s,z}$ over an affine open subscheme $W_z=\Spec(A_z)$ of $Y$ with $z\in W_z$. 

From \cite{B}, Thm. E or \cite{VZ3}, Cor. 3 and the property (ii) we get that for each maximal point $\eta_0\in Y_{\mathbb F_p}$ that belongs to $W_z$, the restriction of $E_{n+s,z}[p^n]$ to $\mathcal O_{Y,\eta_0}$ is (canonically identified with) the truncation of level $n$ of the $\BT$ $D_{\eta_0}$ over $\mathcal O_{Y,\eta_0}$ which extends $D_{\eta}$. Based on Fact \ref{F1} (c) applied to the local rings of $Y$ we get that $D_{Y[\frac{1}{p}]}[p^n]$ and the $E_{n+s,z}[p^n]$'s with $z\in Y_{\mathbb F_p}$ glue together to define a $\BT_n$ $D^+_n$ over $Y$ which extends $D_{\eta}[p^n]$ and whose restriction to each $\mathcal O_{Y,\eta_0}$ is $D_{\eta_0}[p^n]$. 

If $m\in\mathbb N^*$, from Fact \ref{F1} (c) applied to local rings of $Y$ we get that the closed embedding homomorphism $D_{Y[\frac{1}{p}]}[p^n]\rightarrow D_{Y[\frac{1}{p}]}[p^{n+m}]$ extends to a homomorphism $D_n^+\rightarrow D_{n+m}^+$. The fact that the inductive system $D_n^+$ is a $\BT$ over $Y$ which extends $D_{Y[\frac{1}{p}]}$ follows from the fact that its restriction to each $\Spec(\mathcal O_{Y,z})$ is canonically identified with the inductive system $D_z[p^n]$. This ends the proofs of the lemma and of Theorem \ref{T1}. \endproof

\subsection{Proof of Corollary \ref{C1}}\label{S71}

In the situation of Corollary \ref{C1} (a), the condition (i) of Lemma \ref{L8} holds. The condition (ii) of Lemma \ref{L8} also holds: for each affine open cover $(W_{\lambda})_{\lambda\in\Lambda}$ of $Y$ we can take all $N_{\lambda}$'s to be $p-1$, cf. property (ii) of Section \ref{S3}. Thus Lemma \ref{L8} implies that Corollary \ref{C1} (a) holds. In view of the uniqueness (up to unique isomorphism) of the extension from the generic point of $Y$, Corollary \ref{C1} (b) is a particular case of Corollary \ref{C1} (a).\endproof

\section{Complements to Corollary \ref{C1} and Lemma \ref{L8}}\label{S8}

For a topological space $\mathcal Y$, let $\mathcal Y^{\textup{min}}$ be the subspace of maximal points of $\mathcal Y$. We recall that if $\mathcal Y$ is a locally spectral space (e.g., the underlying topological space of a scheme), then $\mathcal Y^{\textup{min}}$ is retrocompact in $\mathcal Y$ if and only if it is proconstructible, and for $\mathcal Y$ spectral if and only if it is quasi-compact, see \cite{ST}, Cors. 2.6 (i) and 2.7.

For the sake of completeness, for extending $\BT$'s from $Y[\frac{1}{p}]$ to $Y$ we have the following variant of Lemma \ref{L8} and Corollary \ref{C1}. 

\begin{lemma}\label{L9}
Let $Y$ be scheme on which $p$ is a non-zero-divisor, which is integrally closed in $Y[\frac{1}{p}]$, and whose local rings at points of residue characteristic $p$ are noetherian.\footnote{Thus in view of \cite{K}, Cor. 7, for each $z\in Y_{\mathbb F_p}$, the local ring $\mathcal O_{Y,z}$ is either a discrete valuation ring or a noetherian ring of both dimension and depth at least $2$. We will only use this characterization and not the integrally closed condition itself.} We assume that the condition (i) of Lemma \ref{L8} holds. Then the following two properties hold:

\medskip
{\bf (a)} If the condition (ii) of Lemma \ref{L8} holds, then each $\BT$ $D_{Y[\frac{1}{p}]}$ over $Y[\frac{1}{p}]$ which extends at each maximal point of $Y_{\mathbb F_p}$, extends to a $\BT$ over $Y$ (uniquely up to unique isomorphism).

\smallskip
{\bf (b)} If $Y^{\textup{min}}_{\mathbb F_p}$ is retrocompact (i.e., the morphism $Y^{\textup{min}}_{\mathbb F_p}\rightarrow Y$ is quasi-compact), then each $\BT$ $D_W$ over an open subscheme $W$ of $Y$ that contains $Y[\frac{1}{p}]$ and $Y^{\textup{min}}_{\mathbb F_p}$ extends to a $\BT$ over $Y$ (uniquely up to unique isomorphism); thus if $Y$ is a faithfully flat regular $\Spec(\mathbb Z_{(p)})$-scheme, then it is $p$-healthy regular.
\end{lemma}

\noindent
{\bf Proof:} Due to the assumptions of the first sentence of the lemma, the functor 
$$(\textup{BT}\;\textup{groups}\;\textup{over}\; Y)\longrightarrow (\textup{BT}\;\textup{groups}\;\textup{over}\; Y[\frac{1}{p}])$$
is fully faithful.  

Thus to prove (a) we can assume that $Y=W_{\lambda}$ is affine. The remaining part of the proof  of (a) is similar to Cases 1 and 2 of the proof of Lemma \ref{L8}, with just one  difference. Referring to the Case 2 of the proof of Lemma \ref{L8}, as we are not assuming $Y[\frac{1}{p}]$ regular, once we obtain $E_{n+s,z}$ over an affine open subscheme $W_z=\Spec( A_z)$, we have to add that by replacing $W_z$ by an affine open subscheme of it we can assume that $E_{n+s,z}[\frac{1}{p}]$ is isomorphic to the restriction of $D_{Y[\frac{1}{p}]}[p^{n+s}]$ to $\Spec(A_z[\frac{1}{p}])$ under an isomorphism which extends the known isomorphism over $\Spec(\mathcal O_{Y,z}[\frac{1}{p}])$.

To prove (b) we can assume that $Y$ is quasi-compact and quasi-separated. By hypotheses, $Y^{\textup{min}}_{\mathbb F_p}$ is quasi-compact. Thus $W$ can be exhausted by quasi-compact open subschemes of it which contain $Y[\frac{1}{p}]\cup Y^{\textup{min}}_{\mathbb F_p}$ and it is enough to prove the extension assertion for each such open subscheme of $W$. Therefore we can assume $W$ is quasi-compact. Let $n\in\mathbb N^*$. Let $z\in Y_{\mathbb F_p}$. From (a) we get that the restriction of $D_W$ to $W\cap\Spec(\mathcal O_{Y,z})$ extends to a $p$-divisible group $D_z$ over $\Spec(\mathcal O_{Y,z})$. A standard limit argument shows that $D_z[p^n]$ spreads out to an extension of $D_W[p^n]_{W_z\cap W}$ to a $\BT$ over $W_z$, where $W_z$ is an affine open neighborhood of $z$ in $Y$. These extensions are unique (up to unique isomorphism) and glue as we have $\iota_*(\mathcal O_W)=\mathcal O_Y$, where $\iota:W\rightarrow Y$ is the open embedding. Thus each $D_W[p^n]$ extends to a $\BT_n$ over $Y$, and by the uniqueness part, these extensions constitute a $\BT$ over $Y$.\endproof

\begin{ex}\label{EX1}
Let $\mathbb Q\subset K_1\subset K_2\subset \cdots$ be a tower of finite field extensions unramified above $p$ such that the union $K_{\infty}:=\cup_{n\in\mathbb N^*} K_n$ has infinitely many places above $p$. Let $O_n$ be the integral closure of $\mathbb Z_{(p)}$ in $K_n$; so $O_{\infty}:=\cup_{n\in\mathbb N^*} O_n$ is the integral closure of $\mathbb Z_{(p)}$ in $K_{\infty}$. Let $\mathfrak{m}\in\Max(O_{\infty})$ be a non-isolated point. For instance, if $K_n:=\mathbb Q(\sqrt{l_1},\ldots,\sqrt{l_n})$, where the $l_i$'s are distinct primes which are squares in $\mathbb Q_p$, then $\Max(O_{\infty})$ is a Cantor set with no isolated point. Let $\pi_n\in O_n$ be a generator of the product of all maximal ideals of $O_n$ different from $\mathfrak{m}\cap O_n$. Let $x$ be an indeterminate and let $A:=\varinjlim O_n[\frac{x}{\pi_n}]$ (a filtered union of subrings of $O_{\infty}[\frac{x}{p}]$). Then the ring $A$ is regular of dimension $2$ and $\Spec(A/pA)^{\textup{min}}$ is not retrocompact, being a topological space homeomorphic to the disjoint union of $\{\mathfrak{m}\}$ and $\Max(O_{\infty})\setminus\{\mathfrak{m}\}$. This example can be modified in a way which allows ramification above $p$ to get examples of faithfully flat $\mathbb Z_{(p)}$-schemes $Y$ which are regular of dimension $2$ and for which $Y^{\textup{min}}$ is not retrocompact and the condition (ii) of Lemma \ref{L8} does not hold. 
\end{ex}

Let $R$ be a regular local ring of mixed characteristic $(0,p)$ and dimension $d\ge 1$. For $e\in\mathbb N^*$ we consider the following condition on $R$:

\medskip
{\bf ($\diamondsuit_e$)} {\it The reduced ring $(R/pR)_{\red}$ is regular and we have $\textup{div}(p)=e\textup{div}(p)_{\red}$ (as divisors on $X=\Spec(R)$).}

\medskip
If condition ($\diamondsuit_e$) holds for $R$, then it also holds for $\widehat{R}$ and for each local ring of $X$ of mixed characteristic $(0,p)$. If $R$ is formally smooth (or only flat with regular special fiber) over a discrete valuation ring $O$ of mixed characteristic $(0,p)$ and absolute ramification index $e$, then condition ($\diamondsuit_e$) holds for $R$. If $e>1$ there are examples in which condition ($\diamondsuit_e$) holds for $\widehat{R}$ but not for $R$ and moreover $R/pR$ is an integral domain. Moreover, we have the following converse which is related to the applicability of Corollary \ref{C1}:

\begin{lemma}\label{L10} 
If condition ($\diamondsuit_e$) holds for $R$ and $p$ does not divide $e$, then $\widehat{R}$ is a formal power series ring over a complete discrete valuation ring $O$ as above.
\end{lemma} 

\noindent
{\bf Proof:} Let $\mathfrak m$ be the maximal ideal of $R$; so $k=R/\mathfrak m$. We consider a Cohen ring $C(k)\subset \widehat{R}$. As $R$ is a unique factorization domain, from the identity $\textup{div}(p)=e\textup{div}(p)_{\red}$ of divisors of $X$ we get that there exist $\pi\in\widehat{R}$ and a unit $u_{\pi}$ of $\widehat{R}$ such that we have $p=\pi^eu_{\pi}$. We write $u_{\pi}=u^{-1}u_1$, where $u\in C(k)$ and $u_1\in 1+\mathfrak m\widehat{R}$. By Hensel's Lemma, $u_1$ has an $e$-th root  $u_1^{\frac{1}{e}}$ in $\widehat{R}$. Thus, by replacing $\pi$ with $\pi u_1^{\frac{1}{e}}$ we can assume that $u_1=1$. As $pu=\pi^e$, the desired discrete valuation subring of $\widehat{R}$ is $O:=C(k)[\pi]=C(k)[x_d]/(x_d^e-pu)$. If $x_1,\ldots,x_{d-1}\in\widehat{R}$ lift a regular system of parameters of $(\widehat{R}/p\widehat{R})_{\red}$, it is easy to see that the natural homomorphism 
$$C(k)[[x_1,\ldots,x_d]]/(x_d^e-pu)=O[[x_1,\ldots,x_{d-1}]]\rightarrow\widehat{R}$$ 
is an isomorphism.\endproof

\bigskip\noindent
{\bf Acknowledgment.}
The second author would like to thank Binghamton and Bielefeld Universities and I. H. E. S., Bures-sur-Yvette for good working conditions; he was partially supported by the NSF grant DMS \#0900967.

\hbox{}
\hbox{Ofer Gabber,\;\;\;Email: gabber@ihes.fr}
\hbox{Address: IH\'ES, Le Bois-Marie, 35, Route de Chartres,} 
\hbox{F-91440 Bures-sur-Yvette, France.}

%\bigskip
\hbox{}
\hbox{Adrian Vasiu,\;\;\;Email: adrian@math.binghamton.edu}
\hbox{Address: Department of Mathematical Sciences, Binghamton University,}
\hbox{Binghamton, P. O. Box 6000, New York 13902-6000, U.S.A.}
\end{document}